\documentclass[11pt]{amsart}
 \newcommand{\ep}{\end{proof}}

\usepackage{latexsym, amssymb, amsthm}
\usepackage[centertags]{amsmath}
\usepackage{m-pictex}

 \newcommand{\sm}{\smallskip}

\newcommand{\ee}{\end{equation}}
\newcommand{\eea}{\end{eqnarray}}
\newcommand{\bean}{\begin{eqnarray*}}
\newcommand{\eean}{\end{eqnarray*}}

 \newif\ifpctex

  \newtheorem{theorem}{Theorem}
  \newtheorem{definition}{Definition}[section]

  \newtheorem{cond}[definition]{Condition}
  \newtheorem{proposition}[definition]{Proposition}
  \newtheorem{lemma}[definition]{Lemma}
  \newtheorem{cor}[definition]{Corollary}
  
  \newcommand{\beCond}[2]{\Rand{\vspace{0,6cm}\tt #1}\begin{cond}[#2]
  \label{#1}} \theoremstyle{definition}
  \newtheorem{remark}[definition]{Remark}
  \newtheorem{example}[definition]{Example}
  \numberwithin{equation}{section}
  \newtheoremstyle{step}{3pt}{0pt}{\itshape}{}{\bf}{}{.5em}{}

\theoremstyle{step} \newtheorem{step}{Step}

\newcommand{\R}{\mathbb{R}} 
\newcommand{\N}{\mathbb{N}} 
\newcommand{\CB}{\mathcal{B}}

\newcommand{\CM}{\mathcal{M}}

 \newcommand{\CZ}{\mathcal{Z}}

\newcommand{\Rand}[1]{\marginpar{#1}} 
\newcommand{\be}[1]{\begin{equation}\label{#1}}
\newcommand{\bew}[1]{\Rand{\vspace{0,6cm}\tt
#1}\begin{equation*}\label{#1}}
\newcommand{\bea}[1]{\Rand{\vspace{0,6cm}\tt
#1}\begin{eqnarray}\label{#1}}

\newcommand{\beL}[2]{\Rand{\vspace{0,6cm}\tt
#1}\begin{lemma}[#2]\label{#1}}
\newcommand{\beD}[2]{\Rand{\vspace{0,6cm}\tt
#1}\begin{definition}[#2]\label{#1}}
\newcommand{\beT}[2]{\Rand{\vspace{0,6cm}\tt
#1}\begin{theorem}[#2]\label{#1}}
\newcommand{\beP}[2]{\Rand{\vspace{0,6cm}\tt
#1}\begin{proposition}[#2]\label{#1}}
\newcommand{\beC}[2]{\Rand{\vspace{0,6cm}\tt
#1}\begin{cor}[#2]\label{#1}}

\begin{document}

\title[Weak convergence of random metric measure spaces]{{\large
  Convergence in distribution of random metric measure spaces}\\[1mm]
  ($\Lambda$-coalescent measure trees)}

\thispagestyle{empty} \author{Andreas Greven}
\address{Andreas Greven\\ Mathematisches Institut\\ University of
  Erlangen\\ Bismarckstr.\ 1$\tfrac 12$\\ D-91054 Erlangen \\ Germany}
  \email{greven@mi.uni-erlangen.de}

\author{Peter Pfaffelhuber}
\address{Peter Pfaffelhuber\\ Zoologisches Institut\\ Ludwig-Maximilian-University Munich
\\ Gro\ss haderner Stra\ss e 2\\ D-82152 Planegg-Martinsried\\
  Germany} \thanks{The research was supported by the
  DFG-Forschergruppe 498 via grant GR 876/13-1,2}
  \email{p.p@lmu.de}

\author{Anita Winter}
\address{Anita Winter\\ Mathematisches Institut\\ University of
  Erlangen\\ Bismarckstr.\ 1$\tfrac 12$\\ D-91054 Erlangen \\ Germany}
  \email{winter@mi.uni-erlangen.de}

\thispagestyle{empty} 



\begin{abstract} We consider the space of complete and separable
  metric spaces which are equipped with a probability measure.
  A notion of convergence is given based on the philosophy
  that a sequence of metric measure spaces converges if and only if
  all finite subspaces sampled from these spaces converge. This
  topology is metrized following Gromov's idea of embedding two metric
  spaces isometrically into a common metric space combined with the
  Prohorov metric between probability measures on a fixed metric
  space. We show that for this topology convergence in distribution
  follows - provided the sequence is tight - from convergence of all
  randomly sampled finite subspaces. We give a characterization of
  tightness based on quantities which are reasonably easy to
  calculate.

  Subspaces of particular
  interest are the space of real trees and of ultra-metric spaces
  equipped with a probability measure.
  As an example we characterize convergence in distribution
  for the (ultra-)metric
  measure spaces given by the random genealogies of the
  $\Lambda$-coalescents.  We show that the $\Lambda$-coalescent
  defines an infinite (random) metric measure space if and only if the
  so-called ``dust-free''-property holds.
\end{abstract}
\maketitle

\thispagestyle{empty}

\section{Introduction and Motivation}
\label{S:intro}
In this paper we study random metric measure spaces which appear frequently in the form of random trees in
probability theory.
Prominent examples are random binary search trees as a special case of
random recursive trees (\cite{DrmotaHwang2005}), ultra-metric
structures in spin-glasses (see, for example,
\cite{BolthausenKistler2006,MPV87}), and coalescent processes in
population genetics (for example, \cite{Hudson1990, Evans2000}). Of
special interest is the continuum random tree, introduced in
\cite{Ald1993}, which is related to several objects, for example,
Galton-Watson trees, spanning trees and Brownian excursions.

Moreover, examples for Markov chains with values in finite trees are
the Aldous-Broder Markov chain which is related to spanning trees
(\cite{Aldous1990}), growing Galton-Watson trees, and
tree-bisection and reconnection which is a method to search through
tree space in phylogenetic
reconstruction (see e.g., \cite{Felsenstein2003}).

Because of the exponential growth of the state space with an
increasing number of vertices tree-valued Markov chains are - even so
easy to construct by standard theory - hard to analyze for their
qualitative properties.  It therefore seems to be reasonable to pass
to a continuum limit and to construct certain limit dynamics and study
them with methods from stochastic analysis.

We will apply this approach in the forthcoming paper
  \cite{GrePfaWin2006} to trees encoding genealogical relationships in
  exchangeable models of populations of constant size. The result will
  be the tree-valued Fleming-Viot dynamics.  For this purpose it is
  necessary to develop systematically the topological properties of
  the state space and the corresponding convergence in distribution.
  The present paper focuses on these topological properties.

As one passes from finite trees to ``infinite'' trees the necessity
arises to equip the tree with a probability measure which allows to
sample typical finite subtrees.
In \cite{Ald1993}, Aldous discusses a notion of convergence in
distribution of a ``consistent'' family of finite random trees towards a
certain limit: the continuum random tree. In order to define
convergence Aldous codes
trees as separable and complete metric spaces satisfying some special
properties for the metric characterizing them as trees which are
embedded into $\ell_1^+$ and equipped with a probability measure.
In this setting finite trees, i.e., trees with finitely many leaves, are always
equipped with the uniform distribution on the set of leaves.
The idea of convergence in distribution of a ``consistent'' family of
finite random trees follows Kolmogorov's theorem which gives
the characterization of convergence of $\R$-indexed stochastic
processes with regular paths.
That is, a sequence has a unique limit provided a tightness condition holds on
path space and assuming that the ``finite-dimensional distributions''
converge.
The analogs of
finite-dimensional distributions are ``subtrees spanned by finitely
many randomly chosen leaves'' and the tightness criterion is built
on characterizations of compact sets in $\ell^+_1$.

Aldous's notion of convergence has been successful for the purpose of
rescaling a given family of trees and showing convergence in
distribution towards a specific limit random tree.  For
example, Aldous shows that suitably rescaled families of critical
finite variance offspring distribution Galton-Watson trees conditioned
to have total population size $N$ converge as $N\to\infty$ to the
Brownian continuum random tree, i.e., the $\R$-tree associated with a
Brownian excursion.  Furthermore, Aldous constructs the genealogical tree
of a resampling population as a metric measure space associated with
the Kingman coalescent, as the limit of $N$-coalescent trees with weight
$1/N$ on each of their leaves.

However, Aldous's ansatz to view trees as closed subsets of
$\ell_1^+$, and thereby using a very particular embedding for the
construction of the topology, seemed not
quite easy and elegant to work with once one wants to construct
tree-valued limit dynamics (see, for example, \cite{EvaPitWin2006},
\cite{EvaWin2006} and \cite{GrePfaWin2006}). More recently, isometry
classes of $\R$-trees, i.e., a particular class of metric spaces, were
introduced, and a means of measuring the distance between two
(isometry classes of) metric spaces were provided based on an
``optimal matching'' of the two spaces yielding the Gromov-Hausdorff
metric (see, for example, Chapter 7 in \cite{BurBurIva01}).

The main emphasis of the present paper is to exploit Aldous's
philosophy of convergence without using Aldous's particular embedding.
That is, we equip the space of separable and complete real
trees which are equipped with a probability measure
with the following topology:
\begin{itemize}
\item  A sequence of trees (equipped with a probability measure)
  converges to a limit
  tree (equipped with a probability measure) if and only if
  all randomly sampled
finite subtrees converge to the corresponding limit subtrees. The
resulting topology is referred to as the {\em Gromov-weak topology}
(compare Definition~\ref{D:00}).
\end{itemize}

Since the construction of the topology works not only for
tree-like metric spaces, but also for the space
of (measure preserving isometry classes of)
metric measure spaces we formulate everything within this framework.
\begin{itemize}
\item
We will see that the Gromov-weak topology on the space of metric measure
spaces is Polish (Theorem~\ref{T:01}).
\end{itemize}
In fact, we metrize the space of metric measure
spaces equipped with the Gromov-weak topology by the
{\em Gromov-Prohorov metric} which combines the two concepts of metrizing
the space of metric spaces and the space of probability measures on a
given metric space in a straightforward way. Moreover, we present a
number of equivalent metrics which might be useful in different contexts.

This then allows to discuss convergence of random variables taking
values in that space.
\begin{itemize}
\item
We next characterize compact sets (Theorem \ref{T:Propprec} combined with Theorem~\ref{T:05}) and
tightness
(Theorem \ref{T:PropTight} combined with Theorem~\ref{T:05}) via quantities which are reasonably  easy to
compute.
\item
We then illustrate
with the example of the $\Lambda$-coalescent tree (Theorem \ref{T:04})
how the tightness characterization  can be applied.
\end{itemize}

We remark that topologies on metric measure spaces are
considered in detail in Section $3\tfrac 12$ of
\cite{Gromov2000}. We are aware that several of our results (in
particular, Theorems~\ref{T:01}, \ref{T:Propprec} and~\ref{T:05})
are stated in \cite{Gromov2000} in a different set-up. While Gromov
focuses on geometric aspects, we provide the tools necessary to do
probability theory on the space of metric measure spaces.
See Remark \ref{Rem:Gromov} for more details on the connection to Gromov's work.

Further related topologies on particular subspaces of isometry
classes of complete and separable metric spaces have already been
considered in \cite{Stu2006} and \cite{EvaWin2006}. Convergence in
these two topologies implies convergence in the {Gromov-weak topology}
but not vice versa.

\subsection*{Outline}
The rest of the paper is organized as follows.
In the next two sections we formulate the main results.
In Section~\ref{S:mms}
we introduce the space of metric measure spaces equipped with
the Gromov-weak topology and characterize their compact
sets. In Section~\ref{S:tigght} we discuss convergence in distribution
and characterize tightness. We then illustrate
the main results introduced so far with the example of the metric measure
tree associated with genealogies generated by
the infinite $\Lambda$-coalescent in
Section~\ref{S:example}.

Sections~\ref{S:GPW} through~\ref{S:equivtop}
are devoted to the proofs of the theorems.  In
Section~\ref{S:GPW} we introduce the Gromov-Prohorov metric as a
candidate for a complete metric which generates the Gromov-weak
topology and show that the generated topology is separable. As a
technical preparation we
collect results on the modulus of mass distribution and the distance
distribution (see Definition~\ref{D:modul}) in Section~\ref{S:massDist}.
In Sections~\ref{S:compact} and~\ref{S:PropTight} we give
characterizations on pre-compactness and tightness for the topology
generated by the Gromov-Prohorov metric.  In Section
\ref{S:equivtop} we prove that the topology generated by the
Gromov-Prohorov metric coincides with the Gromov-weak
topology.

Finally, in
Section~\ref{S:equivMetrics} we provide several other metrics that generate the
Gromov-weak topology.

\section{Metric measure spaces}
\label{S:mms}
As usual, given a topological space $(X,{\mathcal O})$, we denote by
${\mathcal M}_1(X)$ the space of all probability measures on $X$
equipped with the Borel-$\sigma$-algebra ${\mathcal B}(X)$.  Recall
that the support of $\mu$, $\mathrm{supp}(\mu)$, is the smallest
closed set $X_0\subseteq X$ such that $\mu(X\setminus X_0)=0$. The
push forward of $\mu$ under a measurable map $\varphi$ from $X$ into
another metric space $(Z,r_Z)$ is the probability measure
$\varphi_\ast\mu\in{\mathcal M}_1(Z)$ defined by \be{forw}
\varphi_\ast\mu(A) := \mu\big(\varphi^{-1}(A)\big), \ee for all
$A\in{\mathcal B}(Z)$. Weak convergence in $\mathcal M_1(X)$ is denoted by $\Longrightarrow$.
\sm

In the following we focus
on complete and separable metric spaces.
\begin{definition}[Metric measure space]
  A {\em metric measure space} is a complete and separable metric space
  $(X,r)$ which is equipped with a probability measure
  $\mu\in{\mathcal M}_1(X)$. We write $\mathbb{M}$ \label{D:000} for
  the space of measure-preserving isometry classes of complete and
  separable metric measure spaces, where we say that $(X,r,\mu)$ and
  $(X^\prime,r^\prime,\mu^\prime)$ are {measure-preserving} isometric
  if there exists an isometry $\varphi$ between the supports of $\mu$
  on $(X,r)$ and of $\mu'$ on
  $(X^\prime,r')$ such that $\mu'=\varphi_\ast\mu$.
It is clear that the property of being
measure-preserving isometric is an equivalence relation.

We abbreviate $\mathcal X = (X,r,\mu)$ for a whole isometry class of metric
spaces whenever no confusion seems to be possible.\sm
\end{definition}\sm

\begin{remark}
\begin{enumerate}
  \item[{}]
  \item[(i)] Metric measure spaces, or short {\em mm-spaces}, are discussed
    in \cite{Gromov2000} in detail. Therefore they are sometimes also
    referred to as {\em Gromov metric triples} (see, for example,
    \cite{Ver1998}). \label{Rem:07}
  \item[(ii)] We have to be careful to deal with \emph{sets} in the
    sense of the Zermelo-Fraenkel axioms. The reason is that we will
    show in Theorem~\ref{T:01} that $\mathbb{M}$ can be metrized, say by
    $\mathrm{d}$, such that $(\mathbb{M},\mathrm{d})$ is complete and
    separable.  Hence if $\mathbb{P}\in{\mathcal M}_1(\mathbb{M})$ then
    the measure preserving isometry class represented by
    $(\mathbb{M},\mathrm{d},\mathbb{P})$ yields an element in
    $\mathbb{M}$. The way out is to define $\mathbb{M}$ as the space
    of measure preserving isometry classes of those metric spaces
    equipped with a probability measure whose elements are not
    themselves metric spaces. Using this restriction we avoid the
    usual pitfalls which lead to Russell's antinomy.
\hfill$\qed$
\end{enumerate}
\end{remark}\sm

To be in a position to formalize that for a sequence of metric measure
spaces all finite subspaces sampled by the measures sitting on the
corresponding metric spaces converge we next introduce the algebra of
polynomials on $\mathbb{M}$.

\begin{definition}[Polynomials]
  A function $\Phi=\Phi^{n,\phi}:\mathbb{M}\to\R$
\label{D:01}
is called a {\em polynomial (of degree $n$ with respect to the test
function $\phi$)} on
$\mathbb{M}$ if and only if $n\in\N$ is the mimimal number such that there exists a bounded
continuous function $\phi:\,[0,\infty)^{\binom{n}{2}}\to\R$ such that
\be{pp3}
\begin{aligned}
   \Phi\big((X,r,\mu)\big)=
   \int\mu^{\otimes n}(\mathrm{d}(x_1,...,x_n))\,
   \phi\big((r(x_{i},x_{j}))_{1\le i<j\le n}\big),
\end{aligned}
\ee
where $\mu^{\otimes n}$ is the $n$-fold product measure of $\mu$.
Denote by
$\Pi$ the algebra of all polynomials on $\mathbb M$.
\end{definition}\sm

\begin{example} In future work, we are particularly interested in
  tree-like metric spaces, i.e., ultra-metric spaces and $\R$-trees.
  In this setting, functions of the form \label{Exp:00}
  \eqref{pp3} can be, for example, the mean total length or the
  averaged diameter of the sub-tree spanned by $n$ points sampled
  independently according to $\mu$ from the underlying tree.
  \hfill$\qed$\end{example}\sm

The next example  illustrates that one can, of course,  not
separate metric measure spaces by polynomials of degree $2$ only.
\begin{example}
Consider the following two metric measure spaces.
\label{Exp:01}

\beginpicture
\setcoordinatesystem units <1cm,1cm>
\setplotarea x from -1 to 10, y from 1 to 3.7
\plot 1 3 2 2 3 3 2 2 /
\put{$\frac{1}{2}$} [cC] at .6 3
\put{$\frac{1}{2}$} [cC] at 3.4 3
\put{$\bullet$} [cC] at 1 3
\put{$\bullet$} [cC] at 3 3
\put{$\mathcal X$} [cC] at 1.7 1.5

\plot 7 3 8 2 9 3 8 2 8 3 /
\put{$\frac{2-\sqrt{3}}{6}$} [cC] at 6.4 3
\put{$\frac{2+\sqrt{3}}{6}$} [cC] at 9.6 3
\put{$\frac{1}{3}$} [cC] at 8.3 3
\put{$\bullet$} [cC] at 7 3
\put{$\bullet$} [cC] at 9 3
\put{$\bullet$} [cC] at 8 3
\put{$\mathcal Y$} [cC] at 7.7 1.5
\endpicture

Assume that in both spaces the mutual distances between different
points are $1$.
In both cases, the empirical distribution of the distances
  between two points equals $\tfrac{1}{2}\delta_0+\tfrac{1}{2}\delta_{1}$, and hence all polynomials of degree
  $n=2$ agree. But obviously,
$\mathcal X$ and $\mathcal Y$ are not measure preserving
isometric.
\hfill$\qed$\end{example}\sm

The first key observation is that the algebra of polynomials is a rich enough subclass to
determine a metric measure space.

\begin{proposition}[Polynomials separate points] The algebra $\Pi$ of
\label{P:00}
poly\-nomials separates points in $\mathbb{M}$.
\end{proposition}\sm

We need the useful notion of the distance matrix distribution.

\begin{definition}[Distance matrix distribution]
\label{def:distMatDistr}
Let $\mathcal X=(X,r,\mu)\in\mathbb M$ and
the space of infinite (pseudo-)distance matrices
\be{mutdis}
   \mathbb R^{\rm{met}}
 :=
   \big\{(r_{ij})_{1\leq i<j<\infty}:\,
   r_{ij} + r_{jk}\geq r_{ik},\,\forall\,1\leq i<j<k<\infty\big\}.
\ee
Define the map $\iota^{\mathcal X}:\,X^\N\to\mathbb R^{\rm{met}}$ by
\be{e:iota}
   \iota^{\mathcal X}\big(x_1,x_2,...\big):=\big(r(x_i,x_j)\big)_{1\le i<j<\infty},
\ee
and the \emph{distance matrix distribution} of $\mathcal X$ by
\be{distMat}
\nu^{\mathcal X}:=(\iota^{\mathcal X})_\ast\mu^{\otimes \N}.
\ee
\end{definition}
\sm

Note that for $\mathcal X\in\mathbb M$ and $\Phi$ of the form \eqref{pp3}, we have that
\be{ppp3}
\Phi(\mathcal X) = \int \nu^{\mathcal X}\big({\rm d}(r_{ij})_{1\leq i<j}\big) \phi\big((r_{ij})_{1\leq i<j\leq n}\big).
\ee

\begin{proof}[Proof of Proposition \ref{P:00}]
Let $\mathcal X_\ell=(X_\ell,r_\ell,\mu_\ell)\in\mathbb{M}$, $\ell=1,2$,
and assume that $\Phi(\mathcal X_1)=\Phi(\mathcal X_2)$,
for all $\Phi\in\Pi$.
The algebra $\{\phi\in \mathcal C_{\mathrm{b}}(\mathbb
R^{\binom{n}{2}});\,
n\in\mathbb N\}$ is separating in $\mathcal M_1(\mathbb
R^{\text{met}})$ and so $\nu^{\mathcal X_1}=\nu^{\mathcal X_2}$ by \eqref{ppp3}.
Applying Gromov's Reconstruction theorem for mm-spaces (see
Paragraph~$3\tfrac12.5$ in \cite{Gromov2000}),
we find that ${\mathcal X}_1={\mathcal X}_2$.
\end{proof}\sm


We are now in a position to define the Gromov-weak topology.

\begin{definition}[Gromov-weak topology] A sequence $({\mathcal
X}_n)_{n\in\N}$ is said to converge Gromov-weakly to
${\mathcal X}$ in \label{D:00}
$\mathbb{M}$ if and only if
$\Phi({\mathcal X}_n)$ converges to $\Phi({\mathcal X})$ in $\R$, for all
polynomials $\Phi\in\Pi$. We call the corresponding topology ${\mathcal
O}_{\mathbb{M}}$ on
$\mathbb{M}$ the Gromov-weak topology.
\end{definition}\sm

The following result ensures that the state space is suitable to do
probability theory on it.
\begin{theorem}
The space $({\mathbb{M}},{\mathcal
    O}_{\mathrm{\mathbb{M}}})$ is Polish. \label{T:01}
\end{theorem}\sm

In order to obtain later tightness criteria for laws of random
elements in $\mathbb{M}$ we need a
characterization of the compact sets of
$(\mathbb{M},{\mathcal O}_{\mathbb{M}})$.
Informally, a subset of $\mathbb{M}$
will turn out to be pre-compact iff the corresponding sequence of probability
measures put most of their mass on subspaces of a uniformly bounded
diameter, and if the contribution of points which do not carry much
mass in their vicinity is small.

These two criteria lead to the following definitions.

\begin{definition}[Distance distribution and Modulus of mass distribution]
Let $\mathcal X=(X,r,\mu)\in\mathbb{M}$. \label{D:modul}
\begin{enumerate}
\item[(i)] The \emph{distance distribution}, which is an element in
${\mathcal M}_1([0,\infty))$, is given by $w_{\mathcal X}
  := r_\ast \mu^{\otimes 2}$, i.e.,
\be{distInt}
\begin{aligned}
   w_{\mathcal X}(\cdot)
 :=
   \mu^{\otimes 2}\big\{(x,x'): r(x,x') \in \boldsymbol{\cdot}\big\}.
\end{aligned}
\ee
\item[(ii)] For $\delta>0$, define the \emph{modulus of mass distribution} as \be{modul}
   v_\delta(\mathcal X)
 :=
   \inf\Big\{\varepsilon>0:\, \mu\big\{x\in X:\,
   \mu(B_\varepsilon(x))\le\delta\big\}\le\varepsilon\Big\}
\ee
where $B_\varepsilon(x)$ is the open ball with radius $\varepsilon$ and center $x$.
\end{enumerate}
\end{definition}\sm

\begin{remark}\label{Rem:08} Observe that
$w_{\mathcal X}$ and
$v_\delta$ are well-defined because they are constant on
isometry classes of a given metric
measure space.
\end{remark}\sm

The next result characterizes pre-compactness in
$(\mathbb{M},{\mathcal O}_{\mathbb{M}})$.
\begin{theorem}[Characterization of pre-compactness]
A set $\Gamma\subseteq\mathbb{M}$ is
pre-compact in the Gromov-weak topology if and only if the following hold.
\label{T:Propprec}
\begin{itemize}
\item[(i)] The family $\{w_{\mathcal X}:\,\mathcal X\in\Gamma\}$
 is tight.
\item [(ii)] For all $\varepsilon>0$ there exist a
$\delta=\delta(\varepsilon)>0$
such that
\be{e:cpm}
   \sup_{{\mathcal X}\in\Gamma} v_\delta(\mathcal X)
 <
   \varepsilon.
\ee
\end{itemize}
\end{theorem}\sm

\begin{remark}\label{Rem:05}
If $\Gamma=\{\mathcal X_1, \mathcal X_2,...\}$ then
we can replace $\sup$ by $\limsup$ in
\eqref{e:cpm}. \hfill$\qed$
\end{remark}\sm

\begin{example}
In the following we illustrate the two requirements
  for a family in $\mathbb{M}$ to be pre-compact which are given in
  Theorem \ref{T:Propprec} by two counter-examples.
  \label{Exp:star0}
\begin{itemize}
\item[(i)] Consider the isometry classes of the metric measure spaces
  ${\mathcal
    X}_n:=(\{1,2\},r_n(1,2)=n,\mu_n\{1\}=\mu_n\{2\}=\tfrac12)$.  A
  potential limit object would be a metric space with masses
  $\tfrac{1}{2}$ within distance infinity. This clearly does not
  exist.

Indeed,
the family $\{w_{\mathcal X_n}=\tfrac{1}{2}\delta_0+\tfrac{1}{2}\delta_n;\,n\in\N\}$ is
not tight, and hence $\{{\mathcal X}_n;\,n\in\N\}$ is not pre-compact in
$\mathbb{M}$ by Condition (i) of Theorem \ref{T:Propprec}.
\item[(ii)] Consider the isometry classes of the metric measure
spaces $\mathcal X_n=(X_n, r_n, \mu_n)$ given for $n\in\N$ by
\be{eq:Exp:star:1}
  X_n:=\{1,...,2^n\}, \quad r_n(x,y):= \mathbf{1}\{x\neq y\},
  \quad \mu_n:=2^{-n}\sum_{i=1}^{2^n} \delta_i,
\end{equation}
i.e., 
${\mathcal X}_n$ consists of $2^n$
points of mutual distance $1$ and is equipped
with a uniform measure on all points.

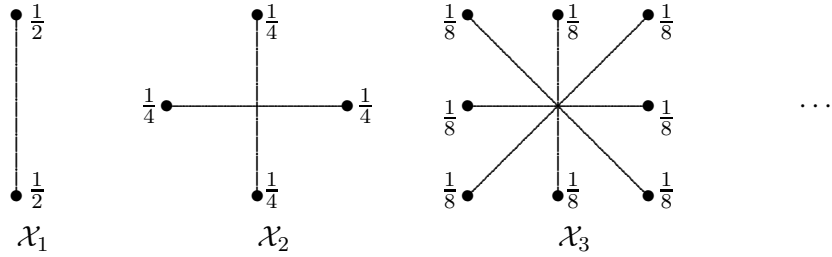
\beginpicture

\setcoordinatesystem units <.8cm,.8cm>
\setplotarea x from -1 to 12, y from -0.5 to 5
\plot 1 1 1 4 /
\put{$\bullet$}[cC] at 1 1
\put{$\bullet$}[cC] at 1 4
\put{$\frac 12$} [lC] at 1.2 1.1
\put{$\frac 12$} [lC] at 1.2 3.9
\put{$\mathcal X_1$} [lC] at 1 0.3

\plot 5 1 5 4 /
\plot 3.5 2.5 6.5 2.5 /
\put{$\bullet$}[cC] at 5 1
\put{$\bullet$}[cC] at 5 4
\put{$\bullet$}[cC] at 3.5 2.5
\put{$\bullet$}[cC] at 6.5 2.5
\put{$\tfrac 14$} [lC] at 5.1 1.1
\put{$\tfrac 14$} [lC] at 5.1 3.9
\put{$\tfrac 14$} [cC] at 3.2 2.5
\put{$\tfrac 14$} [cC] at 6.8 2.5
\put{$\mathcal X_2$} [lC] at 5 0.3

\plot 10 1 10 4 /
\plot 8.5 2.5 11.5 2.5 /
\plot 8.5 1 11.5 4 /
\plot 8.5 4 11.5 1 /
\put{$\bullet$}[cC] at 10 1
\put{$\bullet$}[cC] at 10 4
\put{$\bullet$}[cC] at 8.5 2.5
\put{$\bullet$}[cC] at 11.5 2.5
\put{$\bullet$}[cC] at 11.5 1
\put{$\bullet$}[cC] at 11.5 4
\put{$\bullet$}[cC] at 8.5 1
\put{$\bullet$}[cC] at 8.5 4
\put{$\tfrac 18$} [cC] at 10.25 1.1
\put{$\tfrac 18$} [cC] at 10.25 3.9
\put{$\tfrac 18$} [cC] at 8.2 2.3
\put{$\tfrac 18$} [cC] at 11.8 2.3
\put{$\tfrac 18$} [cC] at 11.8 3.9
\put{$\tfrac 18$} [cC] at 11.8 1.1
\put{$\tfrac 18$} [cC] at 8.2 1.1
\put{$\tfrac 18$} [cC] at 8.2 3.9

\put{$\mathcal X_3$} [lC] at 10 0.3
\put{$\mathbf{\cdots}$} [lC] at 14 2.5
\endpicture

A potential limit object would consist of infinitely many points of mutual distance $1$ with
a uniform measure. Such a space does not exist.

Indeed, notice that for $\delta>0$,
\be{Exp:star:tight2}
    v_\delta(\mathcal X_n)
  =
    \begin{cases} 0, & \delta< 2^{-n},
      \\
    1, & \delta \ge 2^{-n},\end{cases}
\end{equation}
so $\sup_{n\in\mathbb N} v_\delta(\mathcal X_n) = 1$, for all
$\delta>0$. Hence $\{{\mathcal X}_n;\,n\in\N\}$
does not fulfil Condition (ii) of
Theorem~\ref{T:Propprec},
and is therefore not pre-compact.
\hfill$\qed$\end{itemize}
\end{example}\sm

\section{Distributions of random metric measure spaces}
\label{S:tigght}
From Theorem~\ref{T:01} and Definition \ref{D:00} we immediately
conclude the characterization of weak convergence for a sequence of
probability measures on $\mathbb{M}$.
\begin{cor}[Characterization of weak convergence]
  A sequence
  $(\mathbb{P}_n)_{n\in\N}$ in
  ${\mathcal M}_1({\mathbb{M}})$ converges weakly  w.r.t. the Gromov-weak topology
if and only if
\begin{itemize}
\item[(i)] the family $\{\mathbb{P}_n;\,n\in\N\}$
is relatively compact in ${\mathcal M}_1(\mathbb{M})$, and
\item[(ii)] for all polynomials $\Phi\in\Pi$, \label{C:02}
$(\mathbb{P}_n\big[\Phi\big])_{n\in\N}$
converges in $\R$.
\end{itemize}
\end{cor}\sm

\begin{proof} The ``only if'' direction is clear, as polynomials are
  bounded and continuous functions by definition.
  To see the converse, recall from Lemma 3.4.3
  in \cite{EthierKurtz86} that given a relative compact sequence of
  probability measures, each separating family of bounded continuous
  functions is convergence determining.
 \end{proof}\sm

 While Condition (ii) of the characterization of convergence given in
 Corollary~\ref{C:02} can be checked in particular
 examples, we still need a manageable characterization of tightness on
 $\CM_1 (\mathbb{M})$ which we can conclude from
 Theorem~\ref{T:Propprec}. It will be given in terms of the distance
 distribution and the modulus of mass distribution.
\begin{theorem}[Characterization of tightness]
A set $\mathbf{A}\subseteq{\mathcal M}_1(\mathbb{M})$ is
tight if and only if the following holds:
\label{T:PropTight}
\begin{itemize}
\item[(i)] The family $\{\mathbb P[w_\mathcal X]: \mathbb P\in \mathbf A\}$
 is tight in ${\mathcal M}_1(\R)$.
\item [(ii)] For all $\varepsilon>0$ there exist a
$\delta=\delta(\varepsilon)>0$
such that
\be{e:tight}
   \sup_{\mathbb{P}\in\mathbf{A}}\mathbb{P}\big[v_\delta(\mathcal
X)\big]
 <
   \varepsilon.
\ee
\end{itemize}
\end{theorem}\sm

\begin{remark}\label{Rem:06}
\begin{enumerate}
\item[(i)] Using the properties of $v_\delta$ from Lemmata \ref{l:usef} and \ref{l:dvelta} it can be seen that \eqref{e:tight} can be replaced either by
\begin{align}\label{e:tight1a}
\sup_{\mathbb
    P\in\mathbf A} \mathbb P\{v_{\delta}(\mathcal X)\geq
  \varepsilon\}< \varepsilon
\end{align}
or
\begin{align}\label{e:tight1b}
\sup_{\mathbb
    P\in\mathbf A} \mathbb P[\mu\{x: \mu(B_{\varepsilon}(x)) \leq
  \delta\}]< \varepsilon.
\end{align}
\item[(ii)]
If $\mathbf{A}=\{\mathbb{P}_1, \mathbb{P}_2,...\}$ then
we can replace $\sup$ by $\limsup$ in
\eqref{e:tight}, \eqref{e:tight1a} and \eqref{e:tight1b}. \hfill$\qed$
\end{enumerate}
\end{remark}\sm

The usage of Theorem \ref{T:PropTight} will be illustrated
with the example of the {\em $\Lambda$-coalescent measure tree}
constructed in the next section, and with examples of trees corresponding to
spatially structured coalescents
(\cite{GreLimWin2006}) and of evolving coalescents
(\cite{GrePfaWin2006}) in forthcoming work.

\begin{remark} Starting with Theorem~\ref{T:PropTight} one
  characterizes easily \label{Rem:00}
  tightness for the stronger topology given in \cite{Stu2006} based on
certain $L^2$-Wasserstein metrics if one requires in addition to (i)
and (ii) uniform integrability of sampled mutual distance.

Similarly,  with Theorem~\ref{T:PropTight} one
characterizes tightness in the space of measure
  preserving isometry classes of metric spaces equipped with a finite
  measure (rather than a probability measure) if  one requires in
  addition tightness of the family of total masses (compare, also with Remark \ref{Rem:01}(ii)).

\hfill$\qed$\end{remark}\sm

\section{Example: $\Lambda$-coalescent measure trees}
\label{S:example}
In this section we apply the theory of metric measure spaces to a
class of genealogies which arise in population models.  Often such
genealogies are represented by coalescent processes and we focus on
$\Lambda$-coalescents introduced in \cite{Pit1999} (see also
\cite{Sag1999}). The family of $\Lambda$-coalescents appears in the
description of the genealogies of population models with evolution
based on resampling and branching.  Such coalescent processes have
since been the subject of many papers (see, for example,
\cite{MoehleSagitov2001}, \cite{BertoinLeGall2005},
\cite{SevenPeople2005} \cite{LimicSturm2006},
\cite{BereBereSchweinsberg2006}).

In resampling models where the offspring variance of an individual
during a reproduction event is finite, the Kingman coalescent appears
as a special $\Lambda$-coalescent. The fact that general
$\Lambda$-coalescents allow for multiple collisions is reflected in an
infinite variance of the offspring distribution.  Furthermore a
$\Lambda$-coalescent is up to time change dual to the process of
relative frequencies of families of a Galton-Watson process with
possibly infinite variance offspring distribution (compare
\cite{SevenPeople2005}). Our goal here is to decide for which
$\Lambda$-coalescents the genealogies are described by a metric
measure space.

We start with a quick description of $\Lambda$-coalescents.  Recall
that a partition of a set $S$ is a collection $\{A_\lambda\}$ of
pairwise disjoint subsets of $S$, also called {\em blocks}, such that
$S=\cup_\lambda A_\lambda$.  Denote by $\mathbb{S}_\infty$ the
collection of partitions of $\N:=\{1,2,3,...\}$, and for all $n\in\N$,
by $\mathbb{S}_n$ the collection of partitions of $\{1,2,3,...,n\}$.
Each partition $\mathcal{P}\in\mathbb{S}_\infty$ defines an
equivalence relation $\sim_{\mathcal{P}}$ by $i\sim_{\mathcal{P}}j$ if
and only if there exists a partition element $\pi\in\mathcal{P}$ with
$i,j\in\pi$. Write $\rho_n$ for the restriction map from
$\mathbb{S}_\infty$ to $\mathbb{S}_n$. We say that a sequence
$({\mathcal P}_k)_{k\in\N}$ converges in $\mathbb{S}_\infty$ if for
all $n\in\N$, the sequence $(\rho_{n}{\mathcal P}_k)_{k\in\N}$
converges in $\mathbb{S}_n$ equipped with the discrete topology.

We are looking for a strong Markov process
$\xi$ starting in ${\mathcal P}_0\in\mathbb{S}_\infty$
such that for all $n\in\N$, the restricted process
$\xi_n:=\rho_n\circ\xi$ is an $\mathbb{S}_{n
}$-valued Markov chain
which starts in $\rho_n\mathcal{P}_0\in\mathbb{S}_n$,
and given that $\xi_n(t)$
has $b$ blocks, each $k$-tuple of blocks of $\mathbb{S}_n$ is merging
to form a single block  at rate $\lambda_{b,k}$.
Pitman \cite{Pit1999} showed that
such a process exists and is
unique (in law)
if and only if
\be{e:lambda}
   \lambda_{b,k}
 :=
   \int^1_0\Lambda(\mathrm{d}x)\,x^{k-2}(1-x)^{b-k}
\ee
for some non-negative and finite measure $\Lambda$ on the Borel
subsets of $[0,1]$.

Let therefore $\Lambda$ be a non-negative finite measure on $\CB([0,1])$
and $\mathcal{P}\in\mathbb{S}_\infty$. We denote by
$\mathbb{P}^{\Lambda,\mathcal{P}}$
the probability distribution governing $\xi$ with $\xi(0)=\mathcal{P}$ on the
space of cadlag paths with the Skorohod topology.
\begin{example} If we choose
\be{e:triv}
   \mathcal{P}^0
 :=
   \big\{\{1\},\{2\},...\big\},
\ee
$\Lambda=\delta_0$, or
$\Lambda(\mathrm{d}x)=\mathrm{d}x$,
then $\mathbb{P}^{\Lambda,\mathcal{P}^0}$ is the {\em Kingman}
and the {\em Bolthausen-Sznitman} coalescent, respectively.
\hfill$\qed$\end{example}\sm

For each non-negative and finite measure $\Lambda$, all initial
partitions $\mathcal{P}\in\mathbb{S}_{\infty}$ and
$\mathbb{P}^{\Lambda,\mathcal{P}}$-almost all $\xi$, there is a
(random) metric $r^{\xi}$ on $\N$ defined by
\be{e:lambda2}
   r^{\xi}\big(i,j\big)
 :=
   \inf\big\{t\ge 0:\,i\sim_{\xi(t)}j\big\}.
\ee

That is, for a realization $\xi$ of the $\Lambda$ coalescent, $r^{\xi}\big(i,j\big)$ is
the time it needs $i$ and $j$ to coalesce. Notice that
$r^{\xi}$ is an {\em ultra-metric} on $\N$, almost surely, i.e., for all
$i,j,k\in\N$,
\be{e:ultra}
   r^{\xi}(i,j)
 \le
   r^{\xi}(i,k)\vee r^{\xi}(k,j).
\ee

Let $(L^\xi,r^{\xi})$ denote the completion of
$(\N,r^{\xi})$. Clearly, the extension of $r^{\xi}$ to
$L^\xi$ is also an ultra-metric. Recall
that ultra-metric spaces are associated with tree-like
structures.

The main goal of this section is to introduce the {\em
  $\Lambda$-coalescent measure tree} as the metric space
$(L^\xi,r^{\xi})$ equipped with the ``uniform distribution''.
Notice that since the Kingman coalescent is known to ``come down immediately
to finitely many partition elements'' the corresponding metric space
is almost surely compact (\cite{Evans2000}).  Even though there is no
abstract concept of the ``uniform distribution'' on compact spaces,
the reader may find it not surprising that in particular examples one
can easily make sense out of this notion by approximation. We will
see, that for $\Lambda$-coalescents, under an additional assumption on
$\Lambda$, one can extend the uniform distribution to locally compact metric
spaces. Within this class falls, for example, the Bolthausen-Sznitman
coalescent which is known to have infinitely many partition elements
for all times, and whose corresponding metric space is therefore not
compact.

Define $H_n$ to be the map which takes a realization of the
$\mathbb{S}_\infty$-valued coalescent and maps it to (an isometry class of) a
metric measure space as follows:
\be{e:Hn}
   H_n:\,\xi\mapsto
   \Big(L^\xi,r^{\xi},\mu^{\xi}_n:=\tfrac{1}{n}\sum\nolimits_{i=1}^n\delta_i\Big).
\ee

Put then for given ${\mathcal P}_0\in\mathbb{S}_\infty$,
\be{e:Qn}
   \mathbb{Q}^{\Lambda,n}
 :=
   \big(H_n\big)_\ast\mathbb{P}^{\Lambda,\mathcal{P}_0}.
\ee

Next we give the characterization of existence and
uniqueness of the $\Lambda$-coalescent
measure tree.

\begin{theorem}[The $\Lambda$-coalescent measure tree]
The family \label{T:04}
$\{\mathbb{Q}^{\Lambda,n};\,n\in\N\}$
converges in the weak topology with respect to the Gromov-weak topology
if and only if
\be{e:lam}
   \int^1_0 \Lambda(\mathrm{d}x)\,x^{-1}
 =
   \infty.
\ee
\end{theorem}\sm

\begin{remark}[``Dust-free'' property]
  Notice first that Condition (\ref{e:lam}) is equivalent to the total
  coalescence rate of a given $\{i\}\in\mathcal{P}_0$ being infinite
  (compare with the proof of Lemma 25 in \cite{Pit1999}).

By exchangeability and the de Finetti Theorem, the family
$\{\tilde{f}(\pi);\,\pi\in\xi(t)\}$ of frequencies
\be{e:freq}
   \tilde{f}(\pi)
 :=
   \lim_{n\to\infty}\frac{1}{n}\#\big\{j\in\{1,...,n\}:\,j\in\pi\big\}
\ee
exists for $\mathbb{P}^{\Lambda,\mathcal{P}_0}$ almost all $\pi\in\xi(t)$ and
all $t>0$. Define $f:=(f(\pi);\,\pi\in\xi(t))$ to be the ranked
rearrangements of $\{\tilde{f}(\pi);\,\pi\in\xi(t)\}$ meaning that the
entrees of the vector $f$ are non-increasing. Let
$\mathbf{P}^{\Lambda,\mathcal{P}_0}$ denote the probability distribution of
$f$. Call the frequencies $f$ {\em proper} if $\sum_{i\ge
  1}f(\pi_i)=1$.  By Theorem~8 in \cite{Pit1999}, the
$\Lambda$-coalescent has in the limit $n\to\infty$ proper frequencies
if and only if Condition (\ref{e:lam}) holds.

According to {\em Kingman's correspondence} (see, for example,
Theorem~14 in \cite{Pit1999}), the distribution
$\mathbb{P}^{\Lambda,\mathcal{P}_0}$ and $\mathbf{P}^{\Lambda,\mathcal{P}_0}$
determine each other uniquely. For $\mathcal{P}\in\mathbb{S}_\infty$ and
$i\in\N$, let $\mathcal{P}^i:=\{j\in\N:\,i\sim_\mathcal{P} j\}$ denote the partition
element in $\mathcal{P}$ which contains $i$. Then Condition (\ref{e:lam})
holds if and only if for all $t>0$,
\be{zorn7} \mathbb{P}^{\Lambda,\mathcal{P}_0}\big\{
\tilde{f}\big((\xi(t))^1\big)=0\big\} = 0.  \ee

The latter is often referred to as the {\em ``dust''-free} property.
\hfill$\qed$\end{remark}\sm

\begin{proof}[Proof of Theorem~\ref{T:04}]
For {\em existence} we will apply the characterization
  of tightness as given in Theorem~\ref{T:PropTight}, and verify the
  two conditions. \sm

(i) By definition,
for all $n\in\N$, $\mathbb{Q}^{\Lambda,n}[w_{{\mathcal X}}]$ is
exponentially distributed with parameter $\lambda_{2,2}$. Hence the
family $\{\mathbb{Q}^{\Lambda,n}[w_{{\mathcal X}}];\,n\in\N\}$ is tight.
\sm

(ii) Fix $t\in (0,1)$. Then for all $\delta>0$, by the uniform distribution and exchangeability,

\be{e:v1}
\begin{aligned}
   &\mathbb Q^{\Lambda,n}\big[\mu\{x: B_\varepsilon(x)\leq \delta\}]
  \\
 &=
   \mathbb P^{\Lambda,\mathcal P_0}\big[\mu^\xi_n\big\{x\in L^\xi:\, \mu^\xi_n({B}_t(x))\leq\delta
   \big|x=1\big\}]
  \\
 &=
   \mathbb P^{\Lambda,\mathcal P_0}
   \big\{\mu^\xi_n(B_t(1))\leq\delta\big\}.
   \end{aligned}
\ee
By the de Finetti theorem, $\mu^\xi_n(B_t(1))\xrightarrow {n\to\infty}\tilde f\big
   ( (\xi(t))^1\big)$, $\mathbb P^{\Lambda, \mathcal P_0}$-almost
   surely. Hence, dominated convergence yields
\be{e:w2}
\begin{aligned}
   \lim_{\delta\to 0}\lim_{n\to\infty} \mathbb Q^{\Lambda,n}\big[\mu\{x: B_\varepsilon(x)\leq \delta\}]
 &=
   \lim_{\delta\to 0}
   \mathbb P^{\Lambda,\mathcal P_0}\big\{\tilde f((\xi(t))^1)
   \leq\delta\big\}
  \\
 &=
   \mathbb P^{\Lambda,\mathcal P_0}\big\{\tilde f((\xi(t))^1) = 0\big\}.
\end{aligned}
\ee\sm

We have shown that Condition~\eqref{e:lam} is equivalent
to~\eqref{zorn7}, and therefore, using \eqref{e:tight1b}, a limit of $\mathbb Q^{\Lambda, n}$
exists if and only if the ``dust-free''-property holds.\sm

{\em Uniqueness} of the limit points follows from the projective property, i.e.
restricting the observation to a tagged subset of initial individuals
is the same as starting in this restricted initial state.
\end{proof}

\section{A complete metric: The Gromov-Prohorov metric}
\label{S:GPW}
In this section we introduce the Gromov-Prohorov metric $d_{\mathrm{GPr}}$ on
$\mathbb{M}$ and prove that the metric space
$(\mathbb{M},d_{\mathrm{GPr}})$ is complete and separable. In
Section~\ref{S:equivtop} we will see that the Gromov-Prohorov metric generates
the Gromov-weak topology.

Notice that the first naive approach to metrize the Gromov-weak
topology could be to fix a countably dense subset
$\{\Phi_n;\,n\in\N\}$ in the algebra of all polynomials, and to put
for ${\mathcal X},{\mathcal Y}\in\mathbb{M}$, \be{e:naiv}
d_{\mathrm{naive}}\big({\mathcal X},{\mathcal Y}\big) :=
\sum_{n\in\N}2^{-n}\big|\Phi_n({\mathcal X})-\Phi_n({\mathcal
  Y})\big|.  \ee However, such a metric is not complete. Indeed one
can check that the sequence $\{{\mathcal X}_n;\,n\in\N\}$ given in
Example~\ref{Exp:star0}(ii) is a Cauchy sequence w.r.t $d_{\mathrm{naive}}$ which does not
converge.

Recall that metrics on the space of probability measures on a fixed
complete and separable metric space
are well-studied (see, for example, \cite{Rachev1991, GibbsSu2001}).
Some of them, like the Prohorov metric and the Wasserstein metric (on
compact spaces) generate the weak topology. On the other hand the
space of all (isometry classes of compact) metric spaces, not carrying
a measure, is complete and separable once equipped with the
Gromov-Hausdorff metric (see, \cite{EvaPitWin2006}). We recall the
notion of the Prohorov and Gromov-Hausdorff metric below.

Metrics on metric measure spaces should take both components into account and
compare the spaces and the measures simultaneously. This was, for
example, done in \cite{EvaWin2006} and \cite{Stu2006}.  We will follow
along similar lines as in \cite{Stu2006}, but replace
the Wasserstein metric with the Prohorov metric.

~

Recall that the {\em Prohorov metric} between two
probability measures $\mu_1$ and $\mu_2$ on a common metric space $(Z,r_Z)$
is defined by
\be{Proh}
   d_{\mathrm{Pr}}^{(Z,r_Z)}\big(\mu_1,\mu_2\big)
 :=
   \inf\Big\{\varepsilon>0:\,\mu_1(F)\le\mu_2(F^\varepsilon)+
   \varepsilon,\;\forall\,F\text{ closed}\Big\}
\ee
where
\be{eq:Feps}
   F^\varepsilon
 :=
   \big\{z\in Z:\, r_Z(z,z')<\varepsilon,\text{ for some }z'\in F\big\}.
\ee
Sometimes it is easier to work with the equivalent formulation based
on \emph{couplings} of the measures $\mu_1$ and $\mu_2$, i.e.,
measures $\tilde\mu$ on $X\times Y$ with
$\tilde\mu(\boldsymbol{\cdot}\times Y) =\mu_1(\boldsymbol{\cdot})$ and
$\tilde\mu(X\times\boldsymbol{\cdot})=\mu_2(\boldsymbol{\cdot})$.
Notice that the product measure $\mu_1\otimes\mu_2$ is a coupling, and
so the set of all couplings of two measures is not empty. By Theorem
3.1.2 in \cite{EthierKurtz86},
\be{Proh2}
\begin{aligned}
   &d_{\mathrm{Pr}}^{(Z,r_Z)}\big(\mu_1,\mu_2\big)
  \\
 &=
   \inf_{\tilde\mu}\,\inf\Big\{\varepsilon>0:\;
   \tilde\mu\big\{(z,z')\in Z\times Z:\,
   r_Z(z,z')\ge\varepsilon\big\}\le\varepsilon\Big\},
\end{aligned}
\ee
where the infimum is taken over all couplings $\tilde\mu$ of $\mu_1$ and
$\mu_2$. The  metric $d_{\mathrm{Pr}}^{(Z,r_Z)}$ is complete and separable if
$(Z,r_Z)$ is complete and separable (\cite{EthierKurtz86}, Theorem 3.1.7).

~

The Gromov-Hausdorff metric is a metric on the space $\mathbb X_c$ of (isometry classes of) compact metric spaces.
For $(X,r_{X})$ and $(Y,r_{Y})$ in $\mathbb{X}_{\mathrm{c}}$ the \emph{Gromov-Hausdorff metric} is given by
\be{eq:GH1}
   d_{{\mathrm{GH}}}\big((X,r_X),(Y,r_Y)\big)
 :=
   \inf_{(\varphi_{X},\varphi_{Y},Z)}
   d_{\rm{H}}^Z\big(\varphi_X(X),\varphi_Y(Y)\big),
\ee
where the infimum is taken over isometric embeddings $\varphi_X$ and
$\varphi_Y$ from $X$ and $Y$, respectively, into some common metric
space $(Z,r_Z)$, and the
Hausdorff metric $d^{(Z,r_Z)}_{\mathrm{H}}$ for closed subsets of a metric space
$(Z,r_Z)$ is given by
\be{eq:HausDef}
   d_{\rm H}^{(Z,r_Z)}(X,Y)
 :=
   \inf\big\{\varepsilon>0:\,X\subseteq Y^\varepsilon,Y\subseteq X^\varepsilon\big\},
\ee
where $X^\varepsilon$ and $Y^\varepsilon$ are given by
\eqref{eq:Feps}
(compare \cite{Gromov2000, BridsonHaefliger1999, BurBurIva01}).

Sometimes, it is handy to use an equivalent formulation of
the Gromov-Hausdorff metric based on \emph{correspondences}.  Recall
that a \emph{relation} $R$ between two compact metric spaces $(X,r_X)$
and $(Y,r_Y)$ is any subset of $X\times Y$. A relation $R\subseteq
X\times Y$ is called a \emph{correspondence} iff for each $x\in X$
there exists at least one $y\in Y$ such that $(x,y)\in{R}$, and for
each $y'\in Y$ there exists at least one $x'\in X$ such that
$(x',y')\in{R}$. Define the \emph{distortion} of a (non-empty)
relation as
\be{distortion}
   \mathrm{dis}(R)
 :=
   \sup\big\{|r_X(x,x') - r_Y(y,y')|:\, (x,y), (x',y')\in R\big\}.
\end{equation}
Then by Theorem~7.3.25 in \cite{BurBurIva01},
the Gromov-Hausdorff metric can be given in terms of a minimal
distortion of all correspondences, i.e.,
\be{GH}
  d_{{\mathrm{GH}}}\big((X,r_X),(Y,r_Y)\big)
 =
  \frac{1}{2}\inf_{R}{\mathrm{dis}}(R),
\end{equation}
where the infimum is over all correspondences $R$ between $X$ and $Y$.

~

To define a metric between two metric measure spaces $\mathcal
X=(X,r_{X},\mu_X)$ and $\mathcal Y=(Y,r_{Y},\mu_Y)$ in $\mathbb M$, we can neither use the
Prohorov metric nor the Gromov-Hausdorff metric directly. However, we can use the idea due to Gromov
and embed $(X,r_X)$ and $(Y,r_Y)$ isometrically into a common metric
space and measure the distance of the image measures.

\begin{definition}[Gromov-Prohorov metric]
The Gromov-Prohorov distance
between two metric measure spaces $\mathcal X= (X,r_X, \mu_{X})$ and
$\mathcal Y=(Y,r_Y,\mu_{Y})$ in $\mathbb M$ is defined by
\be{dGPr}
   d_{{\mathrm{GPr}}}\big(\mathcal X, \mathcal Y\big)
 :=
   \inf_{(\varphi_{X},\varphi_{Y},Z)} d^{(Z,r_Z)}_{\mathrm{Pr}}
   \big((\varphi_{X})_\ast\mu_{X},(\varphi_{Y})_\ast\mu_{Y}\big),
\ee
where the infimum is taken over all isometric embeddings $\varphi_{X}$ and
$\varphi_Y$ from $X$ and $Y$, respectively, into some common metric space $(Z,r_Z)$.
\end{definition}\sm

\begin{remark}
\begin{enumerate}
\item[(i)]\sloppy To see that the Gromov-Prohorov metric is well-defined we have to check that the right hand side of \eqref{dGPr} does not depend on the element ofthe isometry class of $(X,r_X,\mu_X)$ and $(Y,r_Y,\mu_Y)$. We leave out the straight-forward details.
\item[(ii)] Notice that w.l.o.g.\ the common metric space $(Z,r_Z)$ and the isometric embeddings $\varphi_{X}$ and
$\varphi_Y$ from $X$ and $Y$ can be chosen to be $X\sqcup Y$ and the canonical embeddings $\varphi_X$ and $\varphi_Y$ from $X$ and $Y$ to $X\sqcup Y$, respectively
(compare, for example, Remark~3.3(iii) in \cite{Stu2006}). We can therefore also write
\be{metricdGPr}
   d_{{\mathrm{GPr}}}\big(\mathcal X, \mathcal Y\big)
 :=
   \inf_{r^{X,Y}} d^{(X\sqcup Y,r^{X,Y})}_{\mathrm{Pr}}
   \big((\varphi_X)_\ast\mu_{X},(\varphi_Y)_\ast\mu_{Y}\big),
\ee
where the infimum is here taken over all complete and separable
metrics $r^{X,Y}$ which extend the metrics $r_X$ on $X$ and $r_Y$ on $Y$ to $X\sqcup Y$.
\end{enumerate}
\label{Rem:27}
\end{remark}

\begin{remark}[Gromov's $\underline{\square}_1$-metric]
Even though the material presented in this paper was developed independently
of Gromov's work, some of the most important ideas are already contained in Chapter~3$\tfrac{1}{2}$ in \cite{Gromov2000}.
\label{Rem:Gromov}

More detailed, one can also start with a Polish space $(X,{\mathcal O})$ which is equipped with an probability measure $\mu\in{\mathcal M}_1(X)$ on ${\mathcal B}(X)$, and then introduce a metric $r:X\times X\to\R_+$ as a measurable function
satisfying the metric axioms. Polish measure spaces $(X,\mu)$ can be parameterized by the segment $[0,1)$ where the parametrization refers to a measure preserving map $\varphi:[0,1)\to X$. If $r$ is a metric on $X$ then $r$ can be pulled back to a metric $(\varphi^{-1})_\ast r$ on $[0,1)$ by letting
\be{e:hatr}
   (\varphi^{-1})_\ast r(t,t')
 :=
   r\big(\varphi(t),\varphi(t')\big).
\ee

 Notice that such a measure-preserving parametrization is far from unique and Gromov introduces his $\underline{\square}_1$-distance between $(X,r,\mu)$ and $(X',r',\mu')$ as the infimum of distances $\square_1$ between the two metric spaces
$([0,1),(\varphi^{-1})_\ast r)$ and $([0,1),(\psi^{-1})_\ast r')$ defined as
\be{square}
\begin{aligned}
   &\square_1\big(d,d'\big)
  \\
 &:=
   \sup\big\{\varepsilon>0:\,\exists X_\varepsilon\in{\mathcal B}([0,1)):\,\lambda(X_\varepsilon)\le\varepsilon,\mbox{ s.t.}
  \\
 &\quad\qquad \quad\qquad \quad\qquad
   |d(t_1,t_2)-d'(t_1,t_2)|\le\varepsilon,\;\forall\,t_1,t_2\in X\setminus X_\varepsilon\big\},
\end{aligned}
\ee
where the infimum is taken all possible measure preserving parameterizations
and $\lambda$ denotes the Lebesgue measure.

The interchange of first embedding in a measure preserving way and
then taking the distance between the pulled back metric spaces versus
first embedding isometrically and then taking the distance between the
pushed forward measures explains the similarities between Gromov's
$\varepsilon$-partition lemma (Section~3$\tfrac{1}{2}$.8 in
\cite{Gromov2000}), his union lemma (Section~3$\tfrac{1}{2}$.12 in
\cite{Gromov2000}) and his pre-compactness criterion
(Section~3$\tfrac{1}{2}$.D in \cite{Gromov2000}) on the one hand and
our Lemma~\ref{lZorn}, Lemma~\ref{l:Gpronespace} and
Proposition~\ref{P:02}, respectively, on the other.

We strongly conjecture that the Gromov-weak topology
agrees with the topology generated by Gromov's $\underline{\square}_1$-metric but
a (straightforward) proof is not obvious to us.
\hfill$\qed$
\end{remark}

We first show that the Gromov-Prohorov distance is indeed a metric.
\begin{lemma} $d_{{\mathrm{GPr}}}$ defines a metric on $\mathbb{M}$.
\label{metricL:metric}
\end{lemma}\sm

In the following we refer to the topology generated by the
Gromov-Prohorov metric as the {\em Gromov-Prohorov topology}.
\index{Gromov-Prohorov topology} In
Theorem~\ref{T:05} of Section~\ref{S:equivtop} we will prove that the
Gromov-Prohorov topology and the Gromov-weak topology coincide.
\index{Gromov-weak topology}  \sm

\begin{remark}[Extension of metrics via relations]
\label{rem:met1}
  The proof of the lemma and some of the following results is based on the extension of two metric
  spaces $(X_1,r_{X_1})$ and $(X_2,r_{X_2})$ if a non-empty relation
  $R\subseteq X_1\times X_2$ is known. The result is a metric on
  $X_1\sqcup X_2$ where $\sqcup$ is the disjoint union. Recall the
  distortion of a relation from \eqref{distortion}
Define the metric space $(X_1\sqcup X_2, r_{X_1\sqcup X_2}^R)$ by
letting $r_{{X_1\sqcup X_{2}}}^{R}(x,x'):=r_{X_i}(x,x')$ if $x,x'\in X_i$,
  $i=1,2$ and for $x_1\in X_1$ and $x_2\in X_2$,
\be{def:rZR}
  \begin{aligned}
      &r_{{X_1\sqcup X_{2}}}^{R}(x_1,x_{2})
    \\
   &:=
      \inf\big\{r_{X_1}(x_1,x_1')
      +\tfrac 12\mathrm{dis}(R)+r_{X_{2}}(x_{2},x_{2}'):\,(x_1',x_{2}')\in
      R\big\}.
\end{aligned}
\ee It is then easy to check that $r_{{X_1\sqcup X_{2}}}^{R}$ defines
a (pseudo-)metric on $X_1\sqcup X_{2}$ which extends the metrics on
$X_1$ and $X_{2}$.  In particular, $r^R_{X_1\sqcup
  X_2}(x_1,x_2)=\tfrac 12\text{\rm dis}(R)$, for any pair
$(x_1,x_2)\in R$, and \be{eq:rHZ} d_H^{(X_1\sqcup X_2,r_{X_1\sqcup
    X_2}^R)}(\pi_1 R, \pi_2 R) = \tfrac12\mathrm{dis}(R),
\end{equation}
where $\pi_1$ and $\pi_2$ are the projection operators on $X_1$ and $X_2$, respectively.
\hfill$\qed$
\end{remark}\sm

\begin{proof}[Proof of Lemma \ref{metricL:metric}] {\em Symmetry} 
is obvious and 
  {\em positive definiteness} can be shown by standard
  arguments. 
To see the {\em triangle
    inequality}, let $\varepsilon,\delta> 0$ and ${\mathcal
    X}_i:=(X_i,r_{X_i},\mu_{X_i})\in\mathbb{M}$,
  $i=1,2,3$, be such that $d_{{\mathrm{GPr}}}\big({\mathcal
    X}_1,{\mathcal X}_2\big)<\varepsilon$ and
  $d_{{\mathrm{GPr}}}\big({\mathcal X}_2,{\mathcal X}_3\big)<\delta$.
  Then, by the definition \eqref{dGPr} together with Remark~\ref{Rem:27}(ii), we can find metrics
$r^{1,2}$ and $r^{2,3}$ on $X_1\sqcup X_2$ and $X_2\sqcup X_3$,
  respectively, such that
\be{metric1q1}
   d^{(X_1\sqcup X_2,r^{1,2})}_{\mathrm{Pr}}\big((\tilde\varphi_1)_\ast\mu_{X_1},(\tilde\varphi_2)_\ast\mu_{X_2}\big)<\varepsilon,
\ee
and
\be{metric1q2}
   d^{(X_2\sqcup X_3,r^{2,3})}_{\mathrm{Pr}}
   \big((\tilde\varphi_2')_\ast\mu_{X_2},(\tilde\varphi_3)_\ast\mu_{X_3}\big)<\delta,
\ee
where $\tilde\varphi_1, \tilde\varphi_2$ and $\tilde\varphi_2', \tilde\varphi_3$ are canonical embeddings from $X_1, X_2$ to $X_1\sqcup X_2$ and $X_2, X_3$ to $X_2\sqcup X_3$, respectively. Setting $Z := (X_1\sqcup X_2) \sqcup (X_2\sqcup X_3)$ we define the metric $r_Z^R$ on $Z$ using the relation
\begin{align}
\label{eq:rtriangle}
R := \{(\tilde\varphi_2(x), \tilde\varphi_2'(x)): x\in X_2\} \subseteq (X_1\sqcup X_2)\times (X_2\sqcup X_3)
\end{align}
and Remark \ref{rem:met1}.
Denote the canonical embeddings from $X_1$, the two copies of $X_2$ and $X_3$ to $Z$ by
$\varphi_1, \varphi_2,\varphi_2'$ and $\varphi_3$, respectively. Since $\text{dis}(R)=0$ and
\be{e:trian0}
d_{\text{Pr}}^{(Z,r_Z^R)} \big( (\varphi_2)_\ast \mu_2, (\varphi_2')_\ast \mu_2\big) = 0,
\ee
by the triangle inequality of the Prohorov metric,
\be{metrice:trian}
\begin{aligned}
   d_{{\mathrm{GPr}}}\big({\mathcal X}_1,{\mathcal X}_3\big)
 &\leq
   d_{{\mathrm{Pr}}}^{(Z,r_Z^R)}\big((\varphi_1)_\ast\mu_1,(\varphi_3)_\ast\mu_3\big)
  \\
 &\le
   d_{{\mathrm{Pr}}}^{(Z,r_Z^R)}\big((\varphi_1)_\ast\mu_1,(\varphi_2)_\ast\mu_2\big)+
   d_{{\mathrm{Pr}}}^{(Z,r_Z^R)}\big((\varphi_2)_\ast\mu_2,(\varphi_2')_\ast\mu_2\big) \\ & \qquad \qquad \qquad \qquad \qquad \qquad +
   d_{{\mathrm{Pr}}}^{(Z,r_Z^R)}\big((\varphi_2')_\ast\mu_2,(\varphi_3)_\ast\mu_3\big)
  \\
 &<
   \varepsilon+\delta.
\end{aligned}
\ee
Hence the triangle inequality follows by taking the infimum over all $\varepsilon$ and~$\delta$.
\end{proof}\sm

\begin{proposition} \label{P:07} The metric space is
$(\mathbb{M},d_{\mathrm{GPr}})$ is complete and separable.
\end{proposition}\sm

We prepare the proof with a lemma.
\begin{lemma}\label{Lemm} Fix $(\varepsilon_n)_{n\in\N}$ in
  $(0,1)$.  A sequence
$({\mathcal X}_n:=(X_n,r_{n},\mu_{n}))_{n\in\N}$
in $\mathbb{M}$ satisfies
\be{Lemm1}
   d_{\mathrm{GPr}}\big({\mathcal X}_n,{\mathcal X}_{n+1}\big)
 <
  \varepsilon_n
\ee
if and only if there exist a complete
and separable metric
space $(Z,r_Z)$ and isometric embeddings $\varphi_{1}$, $\varphi_{2}$, ... from
$X_1$, $X_2$, ...,
respectively,
into $(Z,r_Z)$,  such
that
\be{Lemm2}
   d_{\mathrm{Pr}}^{(Z,r_Z)}
   \big((\varphi_{n})_\ast\mu_{n},(\varphi_{{n+1}})_\ast\mu_{{n+1}}
   \big)
 <
  \varepsilon_n.
\ee
\end{lemma}\sm

\begin{proof} The ``if'' direction is clear.  For the ``only if''
  direction, take sequences $({\mathcal
    X}_n:=(X_n,r_{n},\mu_{n}))_{n\in\N}$ and
  $(\varepsilon_n)_{n\in\N}$ which satisfy (\ref{Lemm1}).
  By Remark \ref{Rem:27}, for $Y_n:=X_n\sqcup X_{n+1}$ and all $n\in\N$, there is a metric $r_{Y_n}$ on $Y_n$ such that
  \be{e:ein}
    d_{\rm{Pr}}^{(Y_n,r_{Y_n})} \big((\varphi_{n})_\ast \mu_{n},
  (\varphi_{{n+1}})_\ast \mu_{{n+1}}\big) < \varepsilon_n
  \ee
  where $\varphi_{n}$ and $\varphi_{{n+1}}$ are the canonical embeddings from $X_n$ and $X_{n+1}$ to $Y_n$.
  Put
  \be{e:ein2}
     R_n
   :=
     \big\{(x,x')\in X_n\times X_{n+1}:\,r_{Y_n}(\varphi_{n}(x),
     \varphi_{{n+1}}(x'))< \varepsilon_n\big\}.
  \ee
  Recall from (\ref{Proh2}) that \eqref{e:ein} implies the existence
  of a coupling $\tilde{\mu}_n$ of $(\varphi_{n})_\ast \mu_{n}$ and
  $(\varphi_{{n+1}})_\ast \mu_{{n+1}}$ such that \be{e:ein1}
  \tilde{\mu}_n\big\{(x,x'):\, r_{Y_n}(y,y')<\varepsilon_n\big\} >
  1-\varepsilon_n.  \ee This implies that $R_n$ is not empty and
  \be{eq:rel}
  d_{\text{Pr}}^{(Y_n, r_{Y_n}^{R_n})}
    \big( (\varphi_{n})_\ast\mu_{n}, (\varphi_{{n+1}})_\ast \mu_{{n+1}}\big) \leq \varepsilon_n.
  \ee

Using the metric spaces $(Y_n, r_{Y_n}^{R_n})$ we define recursively metrics $r_{Z_n}$ on
  $Z_n:=\bigsqcup_{k=1}^n X_k$. Starting with $n=1$, we set
  $(Z_1,r_{Z_1}):=({X}_1,r_1)$.  Next, assume we are given a metric
  $r_{Z_n}$ on $Z_n$. Consider the isometric embeddings $\psi_k^n$ from ${X}_k$
  to $Z_n$, for $k=1,...,n$ which arise from the canonical embedding
  of $X_k$ in $Z_n$. Define for all $n\in\N$, \be{e:ein3} \tilde R_n
  := \big\{(z,x_{})\in Z_n\times X_{n+1}: ((\psi_n^n)^{-1}(z),x)\in
  R_n\big\} \ee which defines metrics $r_{Z_{n+1}}^{\tilde R_n}$ on
  $Z_{n+1}$ via \eqref{def:rZR}.

  By this procedure we obtain in the limit a separable metric space
  $(Z':=\bigsqcup_{n=1}^\infty X_n, r_{Z'})$. Denote its completion by
  $(Z, r_{Z})$ and isometric embeddings from $X_n$ to $Z$ which arise by the
  canonical embedding by $\psi_n, n\in\N$. Observe that the
  restriction of $r_{Z}$ to $X_n\sqcup X_{n+1}$ is isometric to
  $(Y_n, r_{Y_n}^{R_n})$ and thus
  \be{eq:rel2}
    d_{\text{Pr}}^{(Z,r_Z)}\big( (\psi_n)_\ast \mu_{X_n}, (\psi_{n+1})_\ast \mu_{X_{n+1}}\big) \leq
  \varepsilon_n
  \ee
  by \eqref{eq:rel}. So the claim follows.
\end{proof}\sm

\begin{proof}[Proof of Proposition~\ref{P:07}]
To get \emph{separability}, we partly
follow the proof of Theorem 3.2.2 in \cite{EthierKurtz86}. Given
${\mathcal X}:=(X,r,\mu)\in\mathbb{M}$ and
$\varepsilon>0$, we can find
${\mathcal
  X}^\varepsilon:=(X,r,\mu^\varepsilon)\in\mathbb{M}$
such that $\mu^\varepsilon$ is a finitely supported atomic measure on $X$ and
$d_{\mathrm{Pr}}(\mu^\varepsilon,\mu)<\varepsilon$. Now
$d_{\mathrm{GPr}}\big({\mathcal X}^\varepsilon,{\mathcal
  X}\big)<\varepsilon$, while $X^\varepsilon$ is just a ``finite
metric space'' and can clearly be approximated arbitrary closely in
the Gromov-Prohorov metric by finite metric spaces with rational
mutual distances and weights.  The set of isometry classes of finite
metric spaces with rational edge-lengths is countable, and so
$(\mathbb{M},d_{{\mathrm{GPr}}})$ is separable.\sm

\sloppy To get {\em completeness}, it suffices to show
that every Cauchy sequence has a convergent subsequence.
Take therefore a Cauchy sequence $({\mathcal X}_n)_{n\in\N}$ in
$(\mathbb{M},d_{\mathrm{GPr}})$ and a subsequence $({\mathcal Y}_n)_{n\in\N}$, $\mathcal Y_n=(Y_n,r_n,\mu_n)$ with $d_{\text{GPr}}(\mathcal Y_n, \mathcal Y_{n+1}) \leq 2^{-n}$.
By Lemma~\ref{Lemm} we can choose  a
complete and separable metric space $({Z},r_{{Z}})$ and,
for each $n\in\N$, an isometric embedding
$\varphi_{n}$ from $Y_n$ into $({Z},r_{{Z}})$
such that $((\varphi_{n})_\ast\mu_{n})_{n\in\N}$
is a Cauchy sequence on ${\mathcal M}_1(Z)$ equipped with the weak
topology. By the completeness of $\mathcal M_1(Z)$,
$((\varphi_{n})_\ast\mu_n)_{n\in\N}$ converges
to some $\bar{\mu}\in\mathcal M_1(Z)$.

Putting the arguments together yields that
with $\mathcal Z:=(Z, r_Z,\bar{\mu})$,
\be{nonum2}
\begin{aligned}
   d_{\rm{GPr}}\big(\mathcal Y_n, \mathcal Z\big)
 \overset{n\to\infty}{\longrightarrow}
   0,
\end{aligned}
\ee
so that $\CZ$ is the desired limit object,
which finishes the proof.
\end{proof}\sm

We conclude this section by another Lemma.

\begin{lemma}
Let $\mathcal X=(X,r,\mu)$, $\mathcal X_1=(X_1,r_{1},\mu_{1})$,
$\mathcal X_2=(X_2,r_{2},\mu_{2}),...$ be in $\mathbb M$.
Then, \label{l:Gpronespace}
\be{sixxt}
   d_{\mathrm{GPr}}\big(\mathcal X_n,\mathcal X\big)
 \xrightarrow{n\to\infty}
   0
\ee
if and only if there exists a complete and separable metric space
$(Z,r_Z)$ and isometric embeddings $\varphi,\varphi_{1},\varphi_{2},...$ from
$X, X_1, X_2$ into $(Z,r_Z)$, respectively,  such that
\be{twee}
   d_{\mathrm{Pr}}^{(Z,r_Z)}\big((\varphi_{n})_\ast\mu_n,\varphi_\ast\mu\big)
 \xrightarrow {n\to\infty}
   0.
\ee
\end{lemma}\sm

\begin{proof} Again the ``if'' direction is clear by definition.
  For the ``only if'' direction, assume that (\ref{sixxt}) holds. To conclude
(\ref{twee}) we can follow the same line of argument as in the proof of Lemma~\ref{Lemm}
but with a metric $r$ extending the metrics $r$, $r_{1}$, $r_2$,... built on correspondences between $X$ and $X_n$
(rather than $X_n$ and $X_{n+1}$). We leave out the details.
\end{proof}

\section{Distance distribution and Modulus of mass distribution}
\label{S:massDist}
In this section we provide results on the distance distribution and on
the modulus of mass distribution. These will be heavily used
in the following sections, where we present metrics which are
equivalent to the Gromov-Prohorov metric and which are very helpful in
proving the characterizations of compactness and tightness in the Gromov-Prohorov topology.

We start by introducing the {\em random distance distribution} of a
given metric measure space.
\begin{definition}[Random distance distribution]
Let $\mathcal X=(X,r,\mu)\in\mathbb M$. For each $x\in X$, define the
map $r_x:X \to [0,\infty)$ by $r_x(x'):=r(x,x')$, and put
\label{def:hatmu}
$\mu^x:=(r_x)_\ast\mu\in \mathcal M_1([0,\infty))$, i.e., $\mu^x$
defines the distribution of distances to the point $x\in
X$. Moreover, define the map $\hat{r}:X\to\mathcal M_1([0,\infty))$
by $\hat{r}(x):=\mu^x$, and let
\be{hatmuuu}
   \hat\mu_{\mathcal X}:=\hat{r}_\ast\mu\in\mathcal M_1(\mathcal
   M_1([0,\infty)))
\ee
be the \emph{random distance distribution} of ${\mathcal X}$.
\end{definition}\sm

Notice first that the random distance
distribution does not characterizes the metric measure space uniquely.
We will illustrate this with an example.

\begin{example}\label{Exp:hatmu}
Consider the
following two metric measure spaces:

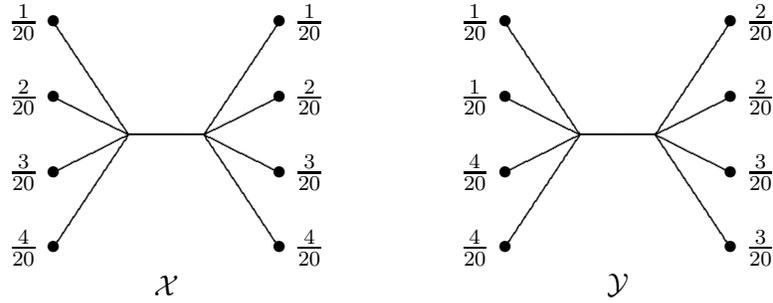
\beginpicture
\setcoordinatesystem units <1cm,1cm>
\setplotarea x from -1 to 10, y from 0 to 5
\plot 1 4 2 2.5 3 2.5 4 4 /
\plot 1 3 2 2.5 3 2.5 4 3 /
\plot 1 2 2 2.5 3 2.5 4 2 /
\plot 1 1 2 2.5 3 2.5 4 1 /
\put{$\frac{1}{20}$} [cC] at .6 4
\put{$\frac{2}{20}$} [cC] at .6 3
\put{$\frac{3}{20}$} [cC] at .6 2
\put{$\frac{4}{20}$} [cC] at .6 1
\put{$\frac{1}{20}$} [cC] at 4.4 4
\put{$\frac{2}{20}$} [cC] at 4.4 3
\put{$\frac{3}{20}$} [cC] at 4.4 2
\put{$\frac{4}{20}$} [cC] at 4.4 1
\put{$\bullet$} [cC] at 1 4
\put{$\bullet$} [cC] at 1 3
\put{$\bullet$} [cC] at 1 2
\put{$\bullet$} [cC] at 1 1
\put{$\bullet$} [cC] at 4 4
\put{$\bullet$} [cC] at 4 3
\put{$\bullet$} [cC] at 4 2
\put{$\bullet$} [cC] at 4 1
\put{$\mathcal X$} [cC] at 2.5 .5

\plot 7 4 8 2.5 9 2.5 10 4 /
\plot 7 3 8 2.5 9 2.5 10 3 /
\plot 7 2 8 2.5 9 2.5 10 2 /
\plot 7 1 8 2.5 9 2.5 10 1 /
\put{$\frac{1}{20}$} [cC] at 6.6 4
\put{$\frac{1}{20}$} [cC] at 6.6 3
\put{$\frac{4}{20}$} [cC] at 6.6 2
\put{$\frac{4}{20}$} [cC] at 6.6 1
\put{$\frac{2}{20}$} [cC] at 10.4 4
\put{$\frac{2}{20}$} [cC] at 10.4 3
\put{$\frac{3}{20}$} [cC] at 10.4 2
\put{$\frac{3}{20}$} [cC] at 10.4 1
\put{$\bullet$} [cC] at 7 4
\put{$\bullet$} [cC] at 7 3
\put{$\bullet$} [cC] at 7 2
\put{$\bullet$} [cC] at 7 1
\put{$\bullet$} [cC] at 10 4
\put{$\bullet$} [cC] at 10 3
\put{$\bullet$} [cC] at 10 2
\put{$\bullet$} [cC] at 10 1
\put{$\mathcal Y$} [cC] at 8.5 .5
\endpicture

That is, both spaces consist of 8 points. The distance between two
points equals the minimal number of edges one has to cross to come
from one point to the other. The measures $\mu_X$ and $\mu_Y$
are given by numbers in the
figure. We find that
\be{tthi}
\begin{aligned}
   \hat\mu_{\mathcal X}
 =
   \hat\mu_{\mathcal Y}
 &=
   \tfrac 1{10}\delta_{\tfrac 1{20}\delta_0+\tfrac{9}{20}\delta_2+\tfrac12\delta_3}+\tfrac {1}{5}\delta_{\tfrac 1{10}
   \delta_0+\tfrac 25\delta_2+\tfrac 1{2}\delta_3}
  \\
 &\quad+
   \tfrac 3{10}\delta_{\tfrac 3{20}\delta_0+\tfrac{7}{20}\delta_2+\tfrac12\delta_3}+\tfrac {2}{5}\delta_{\tfrac 1{5}
   \delta_0+\tfrac 3{10}\delta_2+\tfrac 1{2}\delta_3}.
\end{aligned}
\ee
Hence, the random distance distributions agree. But obviously,
$\mathcal X$ and $\mathcal Y$ are not
measure preserving isometric.
\hfill$\qed$\end{example}\sm

Recall the distance distribution $w_{\boldsymbol{\cdot}}$ and the
modulus of mass
distribution $v_\delta(\boldsymbol{\cdot})$ from
Definition~\ref{D:modul}. Both can be expressed through the random
distance distribution $\hat{\mu}(\boldsymbol{\cdot})$. These
facts follow directly from the definitions, so we omit the
proof.

\begin{lemma}[Reformulation of $w_{\boldsymbol{\cdot}}$ and
  $v_\delta(\boldsymbol{\cdot})$ in terms of $\hat{\mu}(\boldsymbol{\cdot})$]
Let $\mathcal X\in\mathbb M$. \label{l:hatmu1}
\begin{enumerate}
\item[(i)] The distance distribution $w_{\mathcal X}$ satisfies
\be{firs}
   w_{\mathcal X}
 =
   \int_{{\mathcal M}_1([0,\infty))}\hat\mu_{\mathcal X}(\mathrm{d}\nu)\,\nu.
\ee
\item[(ii)] For all $\delta>0$, the modulus of mass
distribution $v_\delta(\mathcal X)$ satisfies
\be{secc}
   v_\delta(\mathcal X)
 =
   \inf\big\{\varepsilon>0: \hat\mu_{\mathcal X}\{\nu\in
   \mathcal M_1([0,\infty)):\,\nu([0,\varepsilon))
   \leq \delta\}\leq\varepsilon\big\}.
\ee
\end{enumerate}
\end{lemma}\sm

The next result will be used frequently.
\begin{lemma} Let ${\mathcal X}=(X,r,\mu)\in\mathbb{M}$ and \label{l:usef}
  $\delta>0$. If $v_\delta({\mathcal X})<\varepsilon$, for some $\varepsilon>0$, then
\be{twew}
   \mu\big\{x\in X:\,\mu(B_\varepsilon(x))\le\delta\big\}
 <
   \varepsilon.
\ee
\end{lemma}\sm

\begin{proof} By definition of $v_\delta(\boldsymbol{\cdot})$, there exists
  $\varepsilon'<\varepsilon$ for which $\mu\big\{x\in X:\,
\mu(B_{\varepsilon'}(x))\le\delta\big\}\le\varepsilon'$.
  Consequently, since $\{x: \mu(B_\varepsilon(x))\leq\delta\} \subseteq
  \{x:\mu(B_{\varepsilon'}(x))\leq\delta\}$,
\be{twone}
  \mu\{x: \mu(B_\varepsilon(x))\leq\delta\} \leq \mu\{x:
  \mu(B_{\varepsilon'}(x))\leq\delta\}
 \le
   \varepsilon'
 <
   \varepsilon,
\ee
and we are done.
\end{proof}\sm

The next result states basic properties of the map $\delta\mapsto v_\delta$.
\begin{lemma}[Properties of $v_\delta(\boldsymbol{\cdot})$]
\label{l:dvelta}
Fix $\mathcal X\in\mathbb{M}$. The map which sends $\delta\ge 0$ to
$v_{\delta}({\mathcal X})$ is non-decreasing, right-continuous and bounded by
$1$. Moreover,
   $v_\delta(\mathcal X)
 \overset{\delta\to 0}{\longrightarrow}
   0$.
\end{lemma}\sm

\begin{proof} The first three properties are trivial. For the forth,
  fix $\varepsilon>0$, and let ${\mathcal
    X}=(X,r,\mu)\in\mathbb{M}$.  Since $X$ is
  complete and separable there exists a compact set
  $K_\varepsilon\subseteq X$ with $\mu(K_\varepsilon)>
  1-\varepsilon$ (see \cite{EthierKurtz86}, Lemma 3.2.1). In
  particular, $K_\varepsilon$ can be covered by finitely many balls
  $A_1,...,A_{N_\varepsilon}$ of radius $\varepsilon/2$ and
  positive $\mu$-mass. Choose $\delta$ such that
\be{p6w}
   0< \delta
 <
   \min\big\{\mu(A_i):\, 1\leq i\leq N_\varepsilon\big\}.
\ee

Then
\be{p6ww}
\begin{aligned}
   \mu\big\{x\in X: \mu(B_\varepsilon(x))>\delta\big\}
 &\ge
   \mu\big(\bigcup\nolimits_{i=1}^{N_\varepsilon}A_i\big)
  \\
 &\ge
   \mu(K_\varepsilon)
  \\
 &>
   1-\varepsilon.
\end{aligned}
\end{equation}

Therefore, by definition, $v_\delta(\mathcal X)\leq\varepsilon$,
and since $\varepsilon$ was chosen arbitrary, the assertion follows.
\end{proof}\sm

The following proposition states
continuity properties of $\hat{\mu}(\boldsymbol{\cdot})$,
$w_{\boldsymbol{\cdot}}$ and $v_\delta(\boldsymbol{\cdot})$.
The reader should have in mind that we finally prove
  with  Theorem~\ref{T:05} in Section~\ref{S:equivtop}
that the Gromov-weak and the
  Gromov-Prohorov topology are the same.
\begin{proposition}[Continuity properties of
  $\hat{\mu}(\boldsymbol{\cdot})$,
  $w_{\boldsymbol{\cdot}}$ and
  $v_\delta(\boldsymbol{\cdot})$]\label{P:lIV2}
\begin{itemize}
\item[{}]
\item[(i)] The map $\mathcal X\mapsto\hat\mu_{\mathcal X}$ is continuous
with respect to the Gromov-weak topology on $\mathbb{M}$ and the weak topology on $\mathcal
M_1(\mathcal M_1([0,\infty)))$.
\item[(ii)] The map $\mathcal X\mapsto\hat\mu_{\mathcal X}$ is continuous
with respect to the Gromov-Prohorov topology on $\mathbb{M}$ and
the weak topology on $\mathcal
M_1(\mathcal M_1([0,\infty)))$.
\item[(iii)]
The map $\mathcal X\mapsto w_{\mathcal X}$
is continuous with respect to both
the Gromov-weak and the Gromov-Prohorov topology on $\mathbb{M}$ and
the weak topology on $\mathcal M_1([0,\infty))$.
\item[(iv)]
Let $\mathcal X$, $\mathcal X_1$, $\mathcal X_2$, ... in $\mathbb M$ such
that $\hat\mu_{\mathcal X_n}\overset{n\to\infty}{\Longrightarrow}
\hat\mu_{\mathcal X}$ and $\delta>0$. Then
\be{sevv}
   \limsup_{n\to\infty}v_{\delta}(\mathcal X_n)
 \leq
   v_\delta(\mathcal X).
\ee
\end{itemize}
\end{proposition}\sm

The proof of Parts~(i) and~(ii) of Proposition~\ref{P:lIV2}
are based on the notion of moment measures.
\begin{definition}[Moment measures of $\hat{\mu}_{\mathcal X}$]\label{D:02}
For $\mathcal X=(X,r,\mu)\in\mathbb M$ and $k\in\N$,
define the  {\em $k^{\mathrm{th}}$ moment measure}
$\hat\mu_{\mathcal X}^k\in\mathcal M_1([0,\infty)^k)$
of $\hat\mu_{\mathcal X}$ by
\be{sixx}
   \hat\mu_{\mathcal X}^k(\mathrm{d}(r_1, ..., r_k))
 :=
   \int\hat\mu_{\mathcal X}(\mathrm{d}\nu)\,
   \nu^{\otimes k}(\mathrm{d}(r_1,...,r_k)).
\ee
\end{definition}\sm

\begin{remark}[Moment measures determine $\hat{\mu}_{\mathcal X}$]
Observe that for all $k\in\N$, \label{Rem:02}
\be{foutt}
\begin{aligned}
   &\hat\mu_{\mathcal X}^k(A_1\times ...\times A_k)
  \\
 &=
   \mu^{\otimes k+1}\big\{(u_0,u_1,...,u_{k}):\,r(u_0, u_1)\in
   A_1,..., r(u_0,u_k)\in A_k\big\}.
\end{aligned}
\ee

By Theorem 16.16 of \cite{Kallenberg2002}, the moment measures
$\hat\mu_{\mathcal X}^k, k=1,2,...$ determine $\hat\mu_{\mathcal X}$
uniquely. Moreover, weak convergence of random measures is equivalent
to convergence of all moment measures. \hfill$\qed$
\end{remark}\sm

\begin{proof}[Proof of Proposition~\ref{P:lIV2}]
(i)
Take $\mathcal X$, $\mathcal X_1$, $\mathcal X_2$, ... in $\mathbb M$
such that
\be{conv}
   \Phi(\mathcal X_n)
 \xrightarrow{n\to\infty}
   \Phi(\mathcal X),
\ee
for all $\Phi\in\Pi$. For $k\in\mathbb N$, consider all
$\phi\in{\mathcal C}_{\mathrm{b}}([0,\infty)^{\binom{k+1}{2}})$ which
depend on $(r_{ij})_{0\le i<j\le k }$ only through
$(r_{0,1}, ..., r_{0,k})$, i.e., there exists $\tilde\phi\in{\mathcal
  C}_{\mathrm{b}}([0,\infty)^{k})$ with $\phi\big((r_{ij})_{0\le
    i<j\le k}\big) =
\tilde\phi\big((r_{0,j})_{1\le j\le k}\big)$. Since for any $\mathcal
Y=(Y,r,\mu)\in\mathbb M$,
\be{fiff}
\begin{aligned}
   &\int\hat\mu^k_{\mathcal Y}(\mathrm{d}(r_1,...,r_k))\,
   \tilde\phi\big(r_1,...,r_k\big)
  \\
 &=
   \int\mu^{\otimes k+1}(\mathrm{d}(u_0,u_1,...,u_{k}))\,
   \tilde\phi\big(r(u_0,u_1),...,r(u_0,u_k)\big)
  \\
 &=
   \int\mu^{\otimes k+1}(\mathrm{d}(u_0,u_1,...,u_{k}))\,
   \phi\big((r(u_i,u_j))_{0\leq i<j\leq k}\big)
\end{aligned}
\ee
it follows from \eqref{conv} that $\hat\mu^k_{\mathcal X_n}\overset{n\to\infty}{\Longrightarrow}\hat\mu_{\mathcal X}^k$ in the topology of
weak convergence. Since $k$ was arbitrary the convergence
$\hat\mu_{\mathcal X_n}\overset{n\to\infty}{\Longrightarrow}\hat\mu_{\mathcal X}$
follows by Remark~\ref{Rem:02}.\sm

(ii) Once more it 
suffices to prove that all moment measures converge. 

Let $\mathcal X=(X,r_X,\mu_X)\in\mathbb{M}$ and
$\varepsilon>0$ be given. Now consider a metric measure space 
$\mathcal Y = (Y,r_Y,\mu_Y) \in\mathbb{M}$ 
with $d_{\rm{GPr}}(\mathcal X, \mathcal Y) < \varepsilon$.

We know that there exists a metric space $(Z,r_Z)$,
isometric embeddings $\varphi_X$ and $\varphi_Y$ of $\mathrm{supp}(\mu_X)$ and $\mathrm{supp}(\mu_Y)$ into
$Z$, respectively, and a coupling $\tilde\mu$ of $(\varphi_X)_\ast \mu_X$ and
$(\varphi_Y)_\ast\mu_Y$ such that
\be{eq:w1}
   \tilde\mu\big\{(z,z'): r_Z(z,z')\ge\varepsilon\big\} 
 \le 
   \varepsilon.
\end{equation}

Given $k\in\mathbb N$, define a coupling $\tilde{\hat\mu}^k$ of $\hat\mu^k_{\mathcal X}$ and $\hat\mu^k_{\mathcal Y}$  
by
\be{yy}
\begin{aligned}
   &\tilde{\hat\mu}^k\big(A_1\times\cdots\times A_k\times
   B_1\cdots\times B_k\big) 
  \\
 &:= 
   \tilde\mu^{\otimes (k+1)}\big\{(z_0,z_0'),...,(z_k,z_k'):\,
   r_Z(z_0,z_i)\in A_i, r_Z(z_0', z_i')\in B_i, i=1,...,k\big\}
\end{aligned}
\end{equation}
for all $A_1\times\cdots\times A_k\times B_1\times\cdots\times B_k\in\mathcal B(\mathbb R_+^{2k})$. Then 
\be{7.10}
\begin{aligned}
   &\tilde{\hat\mu}^k\big\{(r_1,...,r_k,r_1',..., r_k'):\, 
   |r_i-r_i'|\ge 2\varepsilon \text{ for at least one }i\big\}
  \\
 &\leq
   k\cdot\tilde{\hat\mu}^1\big\{(r_1,r_1'):\, |r_1-r_1'|\ge
   2\varepsilon\big\}   
  \\
 &= 
   k\cdot\tilde\mu^{\otimes 2}
   \big\{(z,z'),(\tilde z,\tilde z'):\, |r_Z(z,\tilde z)-r_Z(z',\tilde
   z')|\ge 2\varepsilon\big\} 
  \\
 &\leq 
   k\cdot\tilde\mu^{\otimes 2}\big\{(z,z'),(\tilde z, \tilde z'):\, 
   r_Z(z,z')\ge\varepsilon \text{ or }
   r_Z(\tilde z,\tilde z')\ge\varepsilon\big\} 
  \\
 &\le
   2k\varepsilon,
\end{aligned}
\end{equation}
which implies that $d_{\rm{Pr}}^{\mathbb R_+^k}(\tilde{\hat\mu}^k_{\mathcal
  X}, \tilde{\hat\mu}^k_{\mathcal Y})\le 2k\varepsilon$, and the claim
  follows.\sm

(iii) By Part~(i) of Lemma~\ref{l:hatmu1}, for ${\mathcal
  X}\in\mathbb{M}$, $w_{\mathcal X}$ equals the first moment measure of
  $\hat{\mu}_{\mathcal X}$. The continuity properties
  of ${\mathcal X}\mapsto
  w_{\mathcal X}$ are therefore a direct consequence of (i) and (ii). \sm




(iv) Let $\mathcal X$, $\mathcal X_1$, $\mathcal X_2$, ... in $\mathbb M$ such
that $\hat\mu_{\mathcal X_n}\overset{n\to\infty}{\Longrightarrow}
\hat\mu_{\mathcal X}$ and $\delta>0$. Assume that $\varepsilon>0$ is
such that $\varepsilon>v_\delta(\mathcal X)$. Then by
Lemmata~\ref{l:hatmu1}(ii) and~\ref{l:usef},
\be{eigg}
   \hat\mu_{\mathcal X}\big\{\nu\in\mathcal M_1([0,\infty)):\,
   \nu([0,\varepsilon))\leq\delta\big\}
 <\varepsilon.
\ee

The set $\{\nu\in\mathcal M_1([0,\infty)):\,
\nu([0,\varepsilon))\leq\delta\}$ is closed in $\mathcal
M_1([0,\infty))$. Hence by the Portmanteau Theorem (see, for example, Theorem
3.3.1 in \cite{EthierKurtz86}),
\be{neiii}
\begin{aligned}
   \limsup_{n\to\infty}&\,\hat\mu_{\mathcal X_n}\big\{\nu\in\mathcal
   M_1([0,\infty)):\,\nu([0,\varepsilon))\leq\delta\big\}
  \\
 &\leq
   \hat\mu_{\mathcal X}\big\{\nu\in\mathcal M_1([0,\infty)):\,
   \nu([0,\varepsilon)) \leq\delta\big\}
 <
   \varepsilon.
\end{aligned}
\ee

That is,  we have
$v_\delta(\mathcal X_n)<\varepsilon$, for all but
finitely many $n$, by (\ref{sevv}). Therefore we find that
$\limsup_{n\to\infty} v_\delta(\mathcal X_n)
<\varepsilon$.  This holds for every
$\varepsilon>v_\delta(\mathcal X)$, and
we are done.
\end{proof}\sm

The following estimate will be used in the proofs of
the pre-compactness characterization given in
Proposition~\ref{P:02} and of Part~(i) of Lemma~\ref{LL}.

\begin{lemma}\label{lZorn}
  Let $\delta>0$, $\varepsilon\geq 0$, and $\mathcal
  X=(X,r,\mu)\in\mathbb{M}$. If $v_\delta(\mathcal X)<
  \varepsilon$, then there exists $N\leq \lfloor\frac{1}{\delta}\rfloor$ and
  points $x_1,..., x_N\in X$ such that the following hold.
\begin{itemize}
\item For $i=1,...,N$,
$\mu\big(B_\varepsilon(x_i)\big)>\delta$,
and $\mu\big(\bigcup\limits_{i=1}^N B_{2\varepsilon}(x_i)\big)> 1-\varepsilon$.
\item For all $i,j=1,...,N$ with $i\not=j$,
$r\big(x_i, x_j\big)>\varepsilon$.
\end{itemize}
\end{lemma}\sm

\begin{proof}
  Consider the set $D:=\{x\in X:\,\mu(B_{\varepsilon}(x))>
  \delta\}$. Since $v_\delta(\mathcal X)<\varepsilon$,
  Lemma~\ref{l:usef} implies that $\mu(D)> 1-\varepsilon$.
Take a maximal $2\varepsilon$ separated
  net $\{x_i: i \in I\}\subseteq D$, i.e.,
\be{eq:l12ca}
D\subseteq \bigcup_{i\in I} B_{2\varepsilon}(x_i),
\end{equation}
and for all $i\not =j$,
\be{eq:l12c}
    r(x_i, x_j)>2\varepsilon,
\end{equation}
while adding a further
  point to $D$ would destroy \eqref{eq:l12c}. Such a net exists in every
  metric space (see, for example, in \cite{BurBurIva01},
  p. 278). Since
\be{eq:l13d}
   1
 \geq
   \mu\Big(\bigcup_{i\in I} B_{\varepsilon}(x_i)\Big)
 =
   \sum_{i\in I} \mu\big(B_{\varepsilon}(x_i)\big)
 \geq
   |I|\delta,
\ee
$|I|\leq \lfloor\frac{1}{\delta}\rfloor$ follows.
\end{proof}\sm

\section{Compact sets}
\label{S:compact}
By Prohorov's Theorem, in a complete and separable metric space, a set
of probability measures is relatively compact iff it is tight. This
implies that compact sets in $\mathbb{M}$ play a special role
for convergence results.  In this section we characterize the
(pre-)compact sets in the Gromov-Prohorov topology.

Recall the distance measure $w_{\mathcal X}$ from \eqref{distInt} and
the modulus of mass distribution $v_\delta({\mathcal X})$ from
\eqref{modul}. Denote by $(\mathbb{X}_{\mathrm{c}},d_{\mathrm{GH}})$
the space of all isometry classes of compact metric spaces equipped
with the Gromov-Hausdorff metric (see Section \ref{S:GPW} for basic
definitions).

The following characterizations together with Theorem~\ref{T:05}
stated in Section~\ref{S:equivtop} which
states the equivalence of the Gromov-Prohorov and the Gromov-weak
topology imply the result stated in Theorem~\ref{T:Propprec}.

\begin{proposition}[Pre-compactness characterization] \label{P:02}
Let ${\Gamma}$ be a family in $\mathbb{M}$. The
following four conditions are equivalent.
\begin{itemize}
\item[(a)] The family ${\Gamma}$ is pre-compact in the Gromov-Prohorov
  topology.
\item[(b)] The family $\big\{w({\mathcal X});\,{\mathcal
    X}\in\Gamma\big\}$ is tight, and
\be{eq:unifConvVd}
   \sup_{\mathcal X\in\Gamma}v_\delta({\mathcal X})
 \xrightarrow{\delta \to 0}
   0.
\ee
\item[(c)] For all $\varepsilon>0$ 
there exists $N_\varepsilon\in\N$
such that for all ${\mathcal X}=(X,r,\mu)\in\Gamma$ there is
a subset $X_{\varepsilon, {\mathcal X}}\subseteq X$ with
\begin{itemize}
\item $\mu\big(X_{\varepsilon,{\mathcal X}}\big)\ge 1-\varepsilon$,
\item $X_{\varepsilon, {\mathcal X}}$ can be covered by at most
  $N_\varepsilon$ balls of radius $\varepsilon$, and
\item $X_{\varepsilon,\mathcal X}$ has diameter at most
      $N_\varepsilon$.
\end{itemize}
\item[(d)] For all $\varepsilon>0$
and ${\mathcal X}=(X,r,\mu)\in\Gamma$ there exists
a compact subset $K_{\varepsilon, {\mathcal X}}\subseteq X$ with
\begin{itemize}
\item $\mu\big(K_{\varepsilon,{\mathcal X}}\big)\ge 1-\varepsilon$, and
\item the family $\mathcal{K}_\varepsilon
:=\{K_{\varepsilon,{\mathcal X}};\,
{\mathcal X}\in\Gamma\}$ is pre-compact in
$(\mathbb{X}_{\rm{c}},d_{\mathrm{GH}})$.
\end{itemize}
\end{itemize}
\end{proposition}\sm

\begin{remark} \label{Rem:01}
\begin{itemize}
\item[{}]
\item[(i)]
  In the space of compact metric spaces equipped with a probability
  measure with full support,  Proposition 2.4
  in \cite{EvaWin2006} states that Condition (d) is sufficient for
  pre-compactness.
\item[(ii)] Proposition~\ref{P:02}(b)
  characterizes  tightness for the stronger topology given in \cite{Stu2006} based on
  certain $L^2$-Wasserstein metrics if one requires in addition
  uniform integrability of sampled mutual distance.

  Similarly, (b) characterizes tightness in the space of measure
  preserving isometry classes of metric spaces equipped with a finite
  measure (rather than a probability measure) if  one requires in
  addition tightness of the family of total masses.
\hfill$\qed$
\end{itemize}
\end{remark}\sm

\begin{proof}[Proof of Proposition~\ref{P:02}] As before, we
  abbreviate ${\mathcal X}=(X,r_X,\mu_X)$. We prove four
  implications giving the statement. \sm

  ${(a)\Rightarrow (b).}$ Assume that $\Gamma\in\mathbb{M}$ is
  pre-compact in the Gromov-Prohorov topology.

  To show that $\big\{w({\mathcal X});\,{\mathcal
    X}\in\Gamma\big\}$ is tight, consider a sequence $\mathcal X_1,
\mathcal X_2,...$ in $\Gamma$. Since $\Gamma$ is relatively
compact by assumption, there is a converging subsequence, i.e., we
find $\mathcal X\in\mathbb{M}$ such that
$d_{\rm{GPr}}(\mathcal X_{n_k}, \mathcal X)
\overset{k\to\infty}{\longrightarrow} 0$ along a suitable subsequence
$(n_k)_{k\in\N}$. By Part (iii) of Proposition~\ref{P:lIV2},
$w_{\mathcal X_{n_k}} \overset{k\to\infty}{\Longrightarrow}
w_{\mathcal X}$. As the sequence
was chosen arbitrary it follows that $\big\{w({\mathcal X});\,{\mathcal
    X}\in\Gamma\big\}$ is tight.

The second part of the assertion in (b) is by contradiction. Assume
  that $v_\delta({\mathcal X})$ does not converge to $0$ uniformly in
  ${\mathcal X}\in\Gamma$, as $\delta\to 0$. Then we find an
  $\varepsilon>0$ such that for all $n\in\N$ there exist sequences
  $(\delta_n)_{n\in\N}$ converging to 0 and ${\mathcal X}_n\in\Gamma$ with
\be{eq:P02P1}
   v_{\delta_n}({\mathcal X}_n)
 \ge
   \varepsilon.
\ee

By assumption, there is a subsequence $\{{\mathcal
  X}_{n_k};\,k\in\N\}$, and a metric measure space ${\mathcal
  X}\in\Gamma$ such that $d_{\mathrm{GPr}}\big({\mathcal
  X}_{n_k},{\mathcal X}) \overset{k\to\infty}{\longrightarrow} 0$.  By
Parts~(ii) and~(iv) of Proposition~\ref{P:lIV2},
we find that $\limsup_{k\to\infty} v_{\delta_{n_k}}({\mathcal
  X}_{n_k})=0$
 which contradicts
\eqref{eq:P02P1}. \sm

${(b)\Rightarrow (c).}$
By assumption, for all $\varepsilon>0$ there are $C(\varepsilon)$ with
\be{p9q}
   \sup_{\mathcal X\in\Gamma} w_{\mathcal X}\big([C(\varepsilon),\infty)\big)
 <
   \varepsilon,
\ee
and $\delta(\varepsilon)$ such that
\be{p9qq}
  \sup_{\mathcal X\in\Gamma} v_{\delta(\varepsilon)}({\mathcal X})
 <
  {\varepsilon}.
\ee
Set
\be{Xprime}
   X_{\varepsilon,\mathcal X}'
 :=
   \big\{x\in X:\, \mu_X\big(B_{C(\tfrac{\varepsilon^2}{4})}(x)\big)
   >1-\varepsilon/2\big\}.
\ee

We claim that $\mu_X(X'_{\varepsilon,\mathcal X})>1-\varepsilon/2$.
If this were not the case, there would be $\mathcal X\in\Gamma$ with
\be{eq:IV3d}
\begin{aligned}
   w_{\mathcal X}\big([C(\tfrac14{\varepsilon^2});\infty)\big)
 &=
   \mu_X^{\otimes 2}\big\{(x,x')\in X\times X:\, r_X(x,x')\ge C(\tfrac14{\varepsilon^2})\big\}
  \\
 &\geq
   \mu_X^{\otimes 2}\big\{(x,x'): x\notin X_{\varepsilon,\mathcal X}', x'\notin
   B_{C(\tfrac{\varepsilon^2}{4})}(x)\big\}
  \\
 &\geq
    \frac \varepsilon 2 \mu_X(\complement X_{\varepsilon,\mathcal X}')
  \\
 &\geq
    \frac{\varepsilon^2}{4},
\end{aligned}
\ee
which contradicts (\ref{p9q}).
Furthermore, the diameter of $X'_{\varepsilon,\mathcal X}$ is bounded
by $4C(\tfrac{\varepsilon^2}{4})$.  Indeed, otherwise we would find points $x,x'\in
X'_{\varepsilon,\mathcal X}$ with $B_{C(\tfrac{\varepsilon^2}{4})}(x) \cap B_{C(\tfrac{\varepsilon^2}{4})}(x') =
\emptyset$, which contradicts that
\be{eq:IV3e}
\begin{aligned}
    \mu_X\big(B_{C(\tfrac{\varepsilon^2}{4})}(x)\cap B_{C(\tfrac{\varepsilon^2}{4})}(x')\big)
 &\geq
    1-\mu_X\big(\complement B_{C(\tfrac{\varepsilon^2}{4})}(x)\big)-
    \mu_X\big(\complement B_{C(\tfrac{\varepsilon^2}{4})}(x')\big)
  \\
 &\geq
    1-\varepsilon.
\end{aligned}
\ee\sm

By Lemma~\ref{lZorn}, for all $\mathcal X=(X,r_X,\mu_X)
\in \Gamma$, we can choose points
$x_1,...,x_{N_\varepsilon^{\mathcal X}}\in X$ with
$N_\varepsilon^{\mathcal X} \leq N(\varepsilon)
:=\lfloor \frac{1}{\delta(\varepsilon/2)}\rfloor$,
$r_X(x_i,x_j)>\varepsilon/2$,
$1\le i<j\le N_\varepsilon^{\mathcal X}$,
and
with $\mu_X\big(\bigcup_{i=1}^{N_\varepsilon^{\mathcal X}}
B_{\varepsilon}(x_i)\big) > 1-\varepsilon/2$.

Set
\be{eq:IV3a}
    X_{\varepsilon,\mathcal X}
 :=
    X_{\varepsilon,\mathcal X}'\cap
    \bigcup_{i=1}^{N_\varepsilon^{\mathcal X}}B_{\varepsilon}(x_i).
\ee

Then
$\mu_X(X_{\varepsilon,\mathcal X}) > 1-\varepsilon$. In addition,
$X_{\varepsilon,\mathcal X}$ can be covered by at most $N(\varepsilon)$ balls
of radius $\varepsilon$ and $X'_{\varepsilon, \mathcal X}$ has
diameter at most $4C(\tfrac{\varepsilon^2}{4})$, so the same is true
for $X_{\varepsilon,\mathcal X}$. \sm

${(c)\Rightarrow (d).}$
Fix $\varepsilon>0$, and set $\varepsilon_n
:=\varepsilon 2^{-(n+1)}$, for all $n\in\N$. By assumption we may
choose
for each $n\in\N$, $N_{\varepsilon_n}\in\N$
such that for all ${\mathcal X}\in\Gamma$ there is
a subset $X_{\varepsilon_n, {\mathcal X}}\subseteq X$ of diameter at
most $N_{\varepsilon_n}$ with
$\mu\big(X_{\varepsilon_n,{\mathcal X}}\big)\ge 1-\varepsilon_n$, and such that
$X_{\varepsilon_n, {\mathcal X}}$ can be covered by at most
$N_{\varepsilon_n}$ balls of radius $\varepsilon_n$.
Without loss of generality we may assume
that all $\{X_{\varepsilon_n,\mathcal X};\,n\in\N,{\mathcal
  X}\in\Gamma\}$ are closed.
Otherwise we just take
their closure. For every ${\mathcal X}\in\Gamma$ take compact sets $K_{\varepsilon_n, \mathcal
  X}\subseteq X$ with $\mu_X(K_{\varepsilon_n, \mathcal X})>
1-\varepsilon_n$. Then the set
\be{KepsX}
   K_{\varepsilon,\mathcal X}
 :=
   \bigcap_{n=1}^\infty \big(X_{\varepsilon_n, \mathcal X}
   \cap K_{\varepsilon_n, \mathcal X}\big)
\ee
is compact since it is the intersection of
a compact set with closed sets, and
\be{KepsX2}
   \mu_X(K_{\varepsilon,\mathcal X})
 \geq
   1-\sum_{n=1}^\infty\big(\mu_X(\complement
   X_{\varepsilon_n,\mathcal X})
   +\mu_X(\complement K_{\varepsilon_n,\mathcal X})\big)
 >
   1-\varepsilon.
\ee

Consider
\be{eq:IV3f}
  \begin{aligned}
    \mathcal K_\varepsilon
 :=
    \big\{K_{\varepsilon,\mathcal X};\,\mathcal X\in\Gamma\big\}.
  \end{aligned}
\end{equation}

To show that $\mathcal K_\varepsilon$ is pre-compact we use the pre-compactness criterion given in
Theorem 7.4.15 in \cite{BurBurIva01}, i.e., we have to show that $\mathcal
K_\varepsilon$ is uniformly totally bounded. This means that the
elements of $\mathcal K_\varepsilon$ have bounded diameter and for all
$\varepsilon'>0$ there is a number $N_{\varepsilon'}$ such that all
elements of $\mathcal K_\varepsilon$ can be covered by
$N_{\varepsilon'}$ balls of radius $\varepsilon'$. By definition,
$\mathcal K_{\varepsilon,\mathcal X}\subseteq X_{\varepsilon_1,\mathcal X}$ and so,
$\mathcal K_{\varepsilon,\mathcal X}$ has diameter at most $N_{\varepsilon_1}$. So, take
$\varepsilon'<\varepsilon$ and $n$ large enough for
$\varepsilon_n<\varepsilon'$. Then $X_{\varepsilon_n, \mathcal X}$ as
well as $K_{\varepsilon,\mathcal X}$ can be covered by
$N_{\varepsilon_n}$ balls of radius $\varepsilon'$.
So $\mathcal K_{\varepsilon}$ is pre-compact in
$(\mathbb{X}_{\mathrm{c}},d_{\mathrm{GH}})$. \sm

${(d)\Rightarrow (a).}$ The proof is in two steps. Assume first that
all metric spaces $(X,r_X)$
with $(X,r_X,\mu_X)\in\Gamma$ are compact,
and that the family $\{(X,r_X): (X,r_X,\mu_X) \in\Gamma\}$ is
pre-compact in the Gromov-Hausdorff topology.

Under these assumptions we can choose for every sequence in $\Gamma$ a
subsequence $(\mathcal X_m)_{m\in\N}$, $\mathcal X_m = (X_m, r_{X_m},
\mu_{X_m})$, and a metric space $(X,r_X)$, such that \be{gre5}
d_{\rm{GH}}(X,X_m)\overset{m\to\infty}{\longrightarrow}0.  \ee

By Lemma~\ref{l:met2}, there are a compact metric space
$(Z,r_Z)$ and isometric embeddings $\varphi_X$, $\varphi_{X_1}$, $\varphi_{X_2}$, ...
from $X$, $X_1$, $X_2$, ..., respectively, to $Z$, such that
$d_{\mathrm{H}}\big(\varphi_X(X),
\varphi_{X_m}(X_m)\big)
\overset{m\to\infty}{\longrightarrow}0$.
Since $Z$ is compact, the set $\{(\varphi_{X_m})_\ast\mu_{X_m}: m\in\mathbb N\}$ is pre-compact
in ${\mathcal M}_1(Z)$ equipped with the weak topology.
Therefore $(\varphi_{X_m})_\ast\mu_{X_m}$
has a converging subsequence, and ${(a)}$
follows in this case.

In the second step we consider the general case.  Let
$\varepsilon_n:=2^{-n}$, fix for every ${\mathcal
  X}\in\Gamma$ and every $n\in\N$, $x\in
K_{\varepsilon_n,\mathcal X}$.  Put
\be{gre6}
   \mu_{X,n}(\boldsymbol{\cdot})
 :=
   \mu_X(\boldsymbol{\cdot}\cap K_{\varepsilon_n,\mathcal X})
   +(1-\mu_X(K_{\varepsilon_n,\mathcal
   X}))\delta_x(\boldsymbol{\cdot})
\ee
and let $\mathcal X^n:=(X,r_X,\mu_{\mathcal X,n})$.
By construction, for all $\mathcal X\in\Gamma$,
\be{gre7}
    d_{\rm{GPr}}\big(\mathcal X^n, \mathcal X\big)
 \leq
    \varepsilon_n,
\ee
and $\mu_{X,n}$ is supported by $K_{\varepsilon_n, \mathcal X}$. Hence, $\Gamma_n:=\{\mathcal X^n;\,\mathcal X\in\Gamma\}$
is pre-compact in $\mathbb{X}_{\mathrm{c}}$ equipped with the Gromov-Hausdorff
topology, for all $n\in\N$. We can therefore find a converging
subsequence in $\Gamma_n$, for all $n$, by the first step.

By a diagonal argument we find a subsequence $(\mathcal X_m)_{m\in\N}$ with
$\mathcal X_m=(X_m,r_{X_m},\mu_{X_m})$ such that
$(\mathcal X_m^n)_{m\in\N}$ converges
for every $n\in\N$ to some metric measure space $\mathcal Z_n$.
Pick a subsequence such that for all $n\in\N$ and $m\ge
n$,
\be{gre8}
   d_{\rm{GPr}}\big(\mathcal X_m^n, \mathcal Z_n\big)
 \leq
   \varepsilon_m.
\ee

Then
\be{gre9}
   d_{\rm{GPr}}\big( \mathcal X_m^n, \mathcal X_{m'}^n\big)
 \leq
   2\varepsilon_n,
\ee
for all $m,m'\geq n$.
We conclude that $(\mathcal X_n)_{n\in\N}$ is a Cauchy sequence in
$(\mathbb{M}, d_{\rm{GPr}})$ since $\sum_{n\ge
  1}\varepsilon_n<\infty$.
Indeed,
\be{ttoto}
\begin{aligned}
   &d_{\rm{GPr}}\big(\mathcal X_n, \mathcal X_{n+1}\big)
  \\
 &\leq
   d_{\rm{GPr}}\big(\mathcal X_n, \mathcal X_n^n\big)
  +
  d_{\rm{GPr}}\big(\mathcal X_n^n, \mathcal X_{n+1}^n\big)
  +
  d_{\rm{GPr}}\big(\mathcal X_{n+1}^n, \mathcal X_{n+1}\big)
  \\
 &\leq
  4 \varepsilon_n.
\end{aligned}
\ee

Since $(\mathbb{M}, d_{\rm{GPr}})$ is complete, this sequence
converges and we are done.
\end{proof}\sm

\section{Tightness}
\label{S:PropTight}
In Proposition~\ref{P:02} we have given a characterization for
relative compactness in $\mathbb{M}$ with respect to the
Gromov-Prohorov topology. This characterization extends to the
following tightness characterization in $\mathcal M_1(\mathbb{M})$
which is equivalent to Theorem \ref{T:PropTight}, once we have shown
the equivalence of the Gromov-Prohorov and the Gromov-weak topology in
Theorem~\ref{T:05} in Section~\ref{S:equivMetrics}.

\begin{proposition}[Tightness with respect to the Gromov-Prohorov topology]
  A set $\mathbf{A}\subseteq{\mathcal M}_1(\mathbb{M})$ is tight with
  respect to the Gromov-Prohorov topology on $\mathbb{M}$ if and only
  if for all $\varepsilon>0$ there exist \label{PropTight} $\delta>0$
  and $C>0$ such that
  \be{eq:PropTightPf1} \sup_{\mathbb P\in\mathbf A} \mathbb
  P\big[v_\delta(\mathcal X) + w_{\mathcal X}([C;\infty))\big] <
  \varepsilon.
\end{equation}
\end{proposition}\sm

\begin{proof}[Proof of Proposition~\ref{PropTight}]
For the ``only if'' direction assume that $\mathbf{A}$ is tight and
fix $\varepsilon>0$.  By definition, we find a compact set
$\Gamma_\varepsilon$ in $({\mathbb M},d_{\mathrm{GPr}})$ such that
$\inf_{\mathbb{P}\in\mathbf{A}}\mathbb{P}(\Gamma_\varepsilon)>1-\varepsilon/4$.
Since $\Gamma_\varepsilon$ is compact there are, by part (b) of
Proposition~\ref{P:02}, $\delta=\delta(\varepsilon)>0$ and
$C=C(\varepsilon)>0$ such that $v_\delta(\mathcal X)<\varepsilon/4$
and $w_{\mathcal X}([C,\infty))<\varepsilon/4$, for all $\mathcal X\in
\Gamma_\varepsilon$.  Furthermore
both $v_\delta(\boldsymbol{\cdot})$ and
$w_{\boldsymbol{\cdot}}([C,\infty))$ are bounded above by 1.
Hence for all $\mathbb{P}\in\mathbf{A}$,
\be{pppp}
\begin{aligned}
    &\mathbb{P}\big[v_\delta(\mathcal X)+w_{\mathcal
    X}([C,\infty))\big]
   \\
  &=
    \mathbb{P}\big[v_\delta(\mathcal X)+w_{\mathcal
  X}([C,\infty));\Gamma_\varepsilon\big] +\mathbb{P}\big[v_\delta(\mathcal
  X)+w_{\mathcal X}([C,\infty));\complement \Gamma_\varepsilon\big]
   \\
  &<
     \frac{\varepsilon}{2}+\frac{\varepsilon}{2}
   \\
 &=
   \varepsilon.
\end{aligned}
\end{equation}
Therefore \eqref{eq:PropTightPf1} holds.
\sm

For the ``if'' direction assume \eqref{eq:PropTightPf1} is true and
fix $\varepsilon>0$. For all $n\in\N$, there are $\delta_n>0$ and
$C_n>0$ such that
\be{e:tight1}
\begin{aligned}
  \sup_{\mathbb{P}\in\mathbf{A}}
  \mathbb{P}\big[v_{\delta_n}(\mathcal X)+w_{\mathcal X}([C_n,\infty))\big]
 <
  2^{-2n}\varepsilon^2.
\end{aligned}
\ee

By Tschebychev's inequality, we conclude that for all $n\in\N$,
\be{e:tight2}
    \sup_{\mathbb{P}\in\mathbf{A}}
    \mathbb{P}\big\{\mathcal X:\, v_{\delta_n}(\mathcal X)+
    w_{\mathcal X}([C_n,\infty))> 2^{-n}\varepsilon\big\}
 <
    2^{-n}\varepsilon.
\end{equation}
By the equivalence of (a) and (b) in Proposition~\ref{P:02} the closure of
\be{Gaama}
   \Gamma_\varepsilon
 :=
   \bigcap_{n=1}^\infty\big\{\mathcal X:\,
   v_{\delta_n}(\mathcal X)+w_{\mathcal X}([C_n,\infty))\le
   2^{-n}\varepsilon\big\}
\ee
is compact. We conclude
\be{twel}
\begin{aligned}
  \mathbb{P}\big(\overline{\Gamma_\varepsilon})
 &\ge
  \mathbb{P}\big(\Gamma_\varepsilon)
  \\
 &\ge
   1-\sum_{n=1}^\infty\mathbb{P}\big\{\mathcal X:\,
   v_{\delta_n}(\mathcal X)+w_{\mathcal
 X}([C_n,\infty))>\tfrac{\varepsilon}{2^n}\big\}
  \\
 &>
   1-\varepsilon.
\end{aligned}
\ee
Since $\varepsilon$ was arbitrary, $\mathbf{A}$ is tight.
\end{proof}\sm

\section{Gromov-Prohorov and Gromov-weak topology  coincide}
\label{S:equivtop}
In this section we show that the topologies induced by convergence
of polynomials and
convergence in the Gromov-Prohorov metric
coincide. This implies that  the characterizations of compact subsets of
$\mathbb{M}$  and tight families in ${\mathcal M}_1(\mathbb{M})$ in
Gromov-weak topology stated in Theorems \ref{T:Propprec} and
\ref{T:PropTight} are covered by the corresponding characterizations
with respect to the Gromov-Prohorov topology given in Propositions
\ref{P:02}
and \ref{PropTight}, respectively. Recall the distance matrix distribution from Definition \ref{def:distMatDistr}
\begin{theorem}
Let $\mathcal X, \mathcal X_1,\mathcal X_2,...\in\mathbb{M}$. \label{T:05}
The following are equivalent:
\begin{itemize}
\item[(a)] The Gromov-Prohorov metric converges, i.e.,
\be{e:(b)}
   d_{\mathrm{GPr}}\big(\mathcal X_n,\mathcal X\big)
 \xrightarrow{n\to\infty}
   0.
\ee
\item[(b)] Distance matrix distributions converge, i.e.
\be{e:(0)}
\nu^{\mathcal X_n} \Longrightarrow \nu^{\mathcal X} \text{ as }n\to\infty.
\ee
\item[(c)] All polynomials $\Phi\in\Pi$ converge, i.e.,
\be{e:(a)}
   \Phi(\mathcal X_n)
 \xrightarrow{n\to\infty}
   \Phi(\mathcal X).
\ee
\end{itemize}
\end{theorem}\sm

\begin{proof}
${(a)\Rightarrow (b).}$ Let $\mathcal X=(X,r_X, \mu_X)$, $\mathcal
X_1=(X_1,r_1, \mu_1)$, $\mathcal X_2=(X_2,r_2,\mu_2)$.
By Lemma~\ref{l:Gpronespace} there are a complete
and separable metric space $(Z,r_Z)$ and isometric embeddings $\varphi$,
$\varphi_1$, $\varphi_2$,... from $(X,r_X)$, $(X_1,r_{1})$,
$(X_2,r_{2})$, ..., respectively, to $(Z,r_Z)$ such that
$(\varphi_n)_\ast\mu_n$ converges weakly to $\varphi_\ast\mu_X$ on
$(Z,r_Z)$. Consequently, using \eqref{e:iota},
\be{eq:d1}
\nu^{\mathcal X_n} =
(\iota^{\mathcal X_n})_\ast \mu_n^{\N} =
(\iota^{\mathcal Z})_\ast \big(((\varphi_n)_\ast \mu_n)^{\N}\big) \Longrightarrow
(\iota^{\mathcal Z})_\ast \big((\varphi_\ast \mu)^{\N}\big)
=
(\iota^{\mathcal X})_\ast \mu^{\N} =
\nu^{\mathcal X}
\ee

${(b)\Rightarrow (c).}$
This is a consequence of the two different representation of polynomials from \eqref{pp3} and \eqref{ppp3}.

${(c)\Rightarrow (a).}$ Assume that for all $\Phi\in\Pi$, $\Phi({\mathcal
  X}_n)\xrightarrow{n\to\infty}\Phi({\mathcal X})$.
It is enough to show that the sequence
$({\mathcal X}_n)_{n\in\mathbb N}$ is pre-compact with respect to the
Gromov-Prohorov topology, since by Proposition \ref{P:00}, this would
imply that all limit points coincide and equal ${\mathcal X}$. We
need to check the two conditions guaranteeing pre-compactness given by Part (b)
of Proposition~\ref{P:02}.

By Part~(iii) of Proposition~\ref{P:lIV2}, the map $\mathcal X\mapsto w_{\mathcal X}$ is continuous with
respect to the Gromov-weak topology. Hence, the family $\{w_{\mathcal X_n};\,n\in\N\}$ is
tight.

In addition, by Parts~(i) and (iv) of Proposition~\ref{P:lIV2}, $\limsup_{n\to\infty}
v_{\delta}(\mathcal X_n) \leq v_\delta(\mathcal X)
\xrightarrow{\delta\to 0} 0$. By Remark \ref{Rem:05}, the latter implies
\eqref{eq:unifConvVd}, and we are done. \sm
\end{proof}\sm

\section{Equivalent metrics}
\label{S:equivMetrics}
In Section~\ref{S:GPW} we have seen that $\mathbb{M}$ equipped with
the Gromov-Prohorov metric is separable and complete. In this section we
conclude the paper by presenting
further metrics (not necessarily complete) which are all equivalent to
the Gromov-Prohorov metric and which may be in some situations easier
to work with.

\subsection*{The Eurandom metric}
$\!\!\!\!$\footnote{When we first discussed how to metrize the
Gromov-weak topology the Eurandom metric came up. Since the discussion
took place during a meeting at Eurandom, we decided to name the metric
accordingly.}
Recall from Definition~\ref{D:01} the algebra of polynomials, i.e.,
functions which evaluate distances of finitely
many points sampled from a metric measure space. By
Proposition~\ref{P:00}, polynomials separate
points in $\mathbb{M}$.
Consequently, two metric measure spaces are different if and only if the distributions of
sampled finite subspaces are different.

We therefore define
\be{dGPr2}
\begin{aligned}
   d_{\rm{Eur}}\big(\mathcal X, \mathcal Y\big)
 :=
   \inf_{\tilde\mu}\,\inf\big\{\varepsilon>0:\,
   \tilde{\mu}^{\otimes 2}\{&(x,y),(x',y')\in (X\times Y)^2
   :\;
   \\
    & |r_X(x,x') - r_Y(y,y')|\ge\varepsilon\}<\varepsilon\big\},
\end{aligned}
\end{equation}
where the infimum is over all couplings $\tilde\mu$ of $\mu_X$ and
$\mu_Y$. We will refer to $d_{\rm{Eur}}$ as the {\em Eurandom}
metric.

Not only is $d_{\mathrm{Eur}}$ a metric on $\mathbb{M}$, it also
generates the Gromov-Prohorov topology.

\begin{proposition}[Equivalent metrics]\label{propMet2}
  The distance $d_{\rm{Eur}}$ is a metric on $\mathbb{M}$.
  It is equivalent to $d_{\rm{GPr}}$, i.e., the generated topology
  is the Gromov-weak topology.
\end{proposition}\sm

Before we prove the proposition we give an example to show that the
Eurandom metric is not complete.

\begin{example}[Eurandom metric is not complete]
\label{Exp:star} Let for all $n\in\N$, ${\mathcal
  X}_n:=(X_n,r_n,\mu_n)$ as in Example \ref{Exp:star0}(ii).
For all $n\in\N$,
\be{p3qa}
\begin{aligned}
   d_{\mathrm{Eur}}&\big({\mathcal X}_n,{\mathcal X}_{n+1}\big)
  \\
 &\leq
   \inf\big\{\varepsilon>0:\,(\mu_n\otimes \mu_{n+1})^{\otimes
 2}\{|\mathbf{1}\{x=x'\}-\mathbf{1}\{y=y'\}|\ge\varepsilon\}\le\varepsilon\big\}
  \\
 &\leq
   2^{-(n-1)},
\end{aligned}
\end{equation}
i.e., $(\mathcal X_n)_{n\in\N}$ is a Cauchy sequence for $d_{\rm{Eur}}$ which
does not converge. Hence $(\mathbb{M}, d_{\rm{Eur}})$ is not complete.
The Gromov-Prohorov metric was shown to be complete, and hence the above
sequence is not Cauchy in this metric. Indeed,
\be{p3q}
   d_{\mathrm{GPr}}\big({\mathcal X}_n,{\mathcal X}_{n+1}\big)
 =
   2^{-1}\overset{n\to\infty}{\not\longrightarrow}0.
\end{equation}
\hfill$\qed$\end{example}\sm

To prepare the proof of Proposition~\ref{propMet2},
we provide bounds on the introduced ``distances''.
\begin{lemma}[Equivalence] Let ${\mathcal X},{\mathcal
Y}\in\mathbb{M}$, and $\delta\in(0,\tfrac 12)$. \label{LL}
\begin{itemize}
\item[(i)] If $d_{\mathrm{Eur}}\big({\mathcal X},{\mathcal
Y}\big)<\delta^4$ then $d_{\mathrm{GPr}}\big({\mathcal X},{\mathcal Y}\big)
  <
  12(2v_\delta(\mathcal X)+\delta)$.
\item[(ii)]
\be{eq:eqivMet2}
   d_{\mathrm{Eur}}\big({\mathcal X},{\mathcal Y}\big)
 \leq
   2d_{\mathrm{GPr}}\big({\mathcal X},{\mathcal Y}\big).
\end{equation}
\end{itemize}
\end{lemma}\sm

\begin{proof}
(i) The Gromov-Prohorov metric relies on the
Prohorov metric
of embeddings of $\mu_X$ and $\mu_Y$ in ${\mathcal M}_1(Z)$
in a metric space $(Z,r_Z)$. This
is in contrast to the Eurandom metric which is based on an optimal coupling
of the two measures $\mu_X$ and $\mu_Y$ without referring to a space
of measures over a third
metric space. Since we want to bound the Gromov-Prohorov metric in
terms of the
Eurandom metric
the main goal of the proof is to construct a suitable metric space
$(Z,r_Z)$.

The construction proceeds in three steps.
We start in Step~1 with finding a suitable $\varepsilon$-net $\{x_1,...,x_N\}$
in
$(X,r_X)$, and  show that this net has a suitable
corresponding net $\{y_1,...,y_N\}$ in $(Y,r_Y)$. In Step~2 we then
verify that these nets have the property that
$r_X(x_i, x_j)\approx r_Y(y_i, y_j)$ (where the '$\approx$' is made
precise below) and $\delta$-balls around these nets carry almost all
$\mu_X$- and $\mu_Y$- mass. Finally, in Step~3 we will use these nets
to
define a
metric space $(Z,r_Z)$ containing both  $(X,r_X)$ and $(Y,r_Y)$, and bound
the Prohorov metric of the images of
$\mu_X$ and $\mu_Y$.

\subsubsection*{Step~1 (Construction of suitable $\varepsilon$-nets in
  $X$ and $Y$)} Fix $\delta\in(0,\tfrac{1}{2})$.
Assume that ${\mathcal X},{\mathcal Y}\in\mathbb{M}$ are
  such that $d_{\mathrm{Eur}}\big({\mathcal X},{\mathcal
  Y}\big)<\delta^4$. By definition, we find a coupling $\tilde{\mu}$
  of $\mu_X$ and $\mu_Y$ such that
\be{eq:l11}
\begin{aligned}
  \tilde\mu^{\otimes 2}\big\{(x_1,y_1),(x_2,y_2):\,
  |r_X(x_1,x_2)-r_Y(y_1,y_2)|>2\delta\big\}<\delta^4.
\end{aligned}
\end{equation}

Set $\varepsilon:=4v_\delta(\mathcal X)\geq 0$.  By Lemma~\ref{lZorn},
there are $N\leq\lfloor\frac{1}{\delta}\rfloor$ points $x_1,...,x_N\in
X$ with pairwise distances at least $\varepsilon$,
\be{eq:l12cb}
   \mu\big(B_\varepsilon(x_i)\big)>\delta,
\ee
for all $i=1,...,N$,
and
\be{eq:l12cc}
   \mu\big(\bigcup\nolimits_{i=1}^N B_{\varepsilon}(x_i)\big)
 \geq
   1-\varepsilon.
\end{equation}

Put $D:=\bigcup_{i=1}^N B_{\varepsilon}(x_i)$. We claim that for every
$i=1,...,N$ there is $y_i\in Y$ with
\be{eq:l11aa}
   \tilde\mu\big(B_{\varepsilon}(x_i)\times
   B_{2(\varepsilon+\delta)}(y_i)\big)
 \geq
   (1-\delta^2)\mu_X\big(B_{\varepsilon}(x_i)\big).
\end{equation}

Indeed, assume the
assertion is not true for some $1\leq i\leq N$. Then, for all $y\in
Y$,
\be{eq:l11b}
   \tilde\mu\big(B_\varepsilon(x_i)\times
   \complement B_{2(\varepsilon+\delta)}(y)\big)
 \ge
   \delta^2\mu_X\big(B_\varepsilon(x_i)\big).
\end{equation}
which implies that
\be{eq:l11c}
\begin{aligned}
   \tilde\mu^{\otimes 2}& \{(x',y'), (x'',y''): |r_X(x',x'') -
   r_Y(y',y'')|>2\delta\} \\& \geq \tilde\mu^{\otimes 2}\{(x',y'),
   (x'',y''): x',x''\in B_\varepsilon(x_i), y''\notin
   B_{2(\varepsilon+\delta)}(y')\}
  \\
 &\geq
   \mu_X(B_\varepsilon(x_i))^2\delta^2
 \\
 &>
   \delta^4,
\end{aligned}
\end{equation}
by \eqref{eq:l12ca} and (\ref{eq:l11b}) which contradicts \eqref{eq:l11}.

\subsubsection*{Step~2 (Distortion of $\{x_1,...,x_N\}$ and
$\{y_1,...,y_n\}$)} Assume that $\{x_1,...,x_N\}$ and
$\{y_1,...,y_n\}$ are such that (\ref{eq:l12cb}) through
(\ref{eq:l11aa}) hold.
We claim that then
\be{e:1}
  \big|r_X(x_i,x_j)-r_Y(y_i,y_j)\big|
 \le
  6(\varepsilon+\delta),
\end{equation}
for all $i,j=1,...,N$. Assume that \eqref{e:1} is not true for some
pair $(i,j)$.  Then for all $x'\in B_\varepsilon(x_i)$, $x''\in
B_\varepsilon(x_j)$, $y'\in B_{2(\varepsilon+\delta)}(y_i)$, and $y''\in
B_{2(\varepsilon+\delta)}(y_j)$,
\be{thenfor}
    \big|r_X(x',x'') - r_Y(y',y'')\big|
 >
   6(\varepsilon+\delta)-2\varepsilon-4(\varepsilon+\delta)
 =
   2\delta.
\end{equation}

Then
\be{eq:l12d}
\begin{aligned}
    \tilde\mu^{\otimes 2}&\big\{(x',y'), (x'',y''): |r_X(x',x'') -
    r_Y(y',y'')|>2\delta\big\}
  \\
 &\geq
    \tilde\mu^{\otimes 2}\big\{(x',y'), (x'',y''):
        \\& \qquad \quad
       x'\in B_\varepsilon(x_i),x''\in B_\varepsilon(x_j),
       y'\in B_{2(\varepsilon+\delta)}(y_i),
       y''\in B_{2(\varepsilon+\delta)}(y_j)\}
  \\
 &=
    \tilde\mu\big(B_\varepsilon(x_i)\times
    B_{2(\varepsilon+\delta)}(y_i)\big)
    \tilde\mu\big(B_\varepsilon(x_j)\times B_{2(\varepsilon+\delta)}
    (y_j)\big)
  \\
 &>
     \delta^2(1-\delta)^2
   \\
 &>
     \delta^4,
\end{aligned}
\end{equation}
where we used \eqref{eq:l11aa}, \eqref{eq:l12ca} and $\delta<\tfrac
12$. Since (\ref{eq:l12d}) contradicts \eqref{eq:l11}, we are done.

\subsubsection*{Step~3 (Definition of a suitable metric space $(Z,r_Z)$)}
Define the relation $R:=\{(x_i,y_i): i=1,...,N\}$ between $X$ and
$Y$ and consider the metric space $(Z,r_Z)$ defined by $Z:=X\sqcup Y$
and $r_Z:=r_{X\sqcup Y}^R$, given as in Remark ~\ref{rem:met1}. Choose isometric embeddings
$\varphi_X$ and $\varphi_Y$ from $(X,r_X)$ and $(Y,r_Y)$, respectively, into $(Z,r_Z)$. As
$\mathrm{dis}(R)\leq 6(\varepsilon+\delta)$ (see (\ref{distortion})
for definition), by Remark~\ref{rem:met1},
$r_Z(\varphi_X(x_i),\varphi_Y(y_i))\leq
3(\varepsilon+\delta)$, for all $i=1,...,N$.

If $x\in X$ and $y\in Y$ are such that
$r_Z(\varphi_X(x),\varphi_Y(y)) \geq 6(\varepsilon + \delta)$ and
$r_X(x,x_i)<\varepsilon$ then
\be{eq:l15a}
\begin{aligned}
   r_Y(y,y_i)
 &\geq
   r_Z(\varphi_X(x),\varphi_Y(y)) - r_X(x,x_i) -
   r_Z(\varphi_X(x_i), \varphi_Y(y_i))
  \\
 &\geq
   6(\varepsilon+\delta)-\varepsilon-3(\varepsilon+\delta)
  \\
 &\ge
   2(\varepsilon+\delta)
\end{aligned}
\end{equation}
and so for all $x\in B_{\varepsilon}(x_i)$,
\be{l16}
   \big\{y\in Y:\, r_Z(\varphi_X(x),\varphi_Y(y))\ge6(\varepsilon+\delta)
   \big\}
 \subseteq
   \complement B_{2(\varepsilon+\delta)}(y_i).
\end{equation}

Let $\tilde{\mu}$ be the probability measure on $Z\times Z$ defined by
$\hat{\mu}(A\times
B):=\tilde{\mu}(\varphi^{-1}_{X}(A)\times\varphi^{-1}_{Y}(B))$, for
all $A,B\in{\mathcal B}(Z)$.  Therefore, by \eqref{e:1}, (\ref{l16}),
(\ref{eq:l11aa}) and as $N\le\lfloor 1/\delta\rfloor$,
\be{l13}
\begin{aligned}
   \hat\mu\{&(z,z'): r_Z(z,z')\geq 6(\varepsilon+\delta)\}
  \\
 &\leq
   \hat\mu\big( \varphi_X(\complement D)\times \varphi_Y(Y))+\hat\mu\Big(\bigcup\limits_{i=1}^N
   B_{\varepsilon}(\varphi_{X}(x_i))\times\complement
   B_{2(\varepsilon+\delta)}(\varphi_{Y}(y_i))\Big)
  \\
 &\leq
   \varepsilon+\sum_{i=1}^N\mu_X(B_{\varepsilon} (x_i))\delta^2
  \\ &
 \leq
   \varepsilon+\delta.
\end{aligned}
\end{equation}

Hence, using
(\ref{Proh2}) and $\varepsilon=4v_\delta(\mathcal X)$,
\be{l14}
   d_{\mathrm{Pr}}^{(Z,r_Z)}\big((\varphi_{X})_\ast\mu_X,
   (\varphi_{Y})_\ast\mu_Y\big)
 \le
   6\big(4v_\delta({\mathcal X})+2\delta\big),
\ee
and so $d_{\mathrm{GPr}}\big({\mathcal X},{\mathcal Y}\big)\le
12\big(2v_\delta({\mathcal X})+\delta\big)$, as claimed.
\sm

(ii) Assume that
  $d_{\mathrm{GPr}}\big({\mathcal X},{\mathcal Y}\big)<\delta$. Then,
  by definition, there exists a metric space $(Z,r_Z)$, isometric embeddings
  $\varphi_{X}$ and $\varphi_{Y}$ between $\mathrm{supp}(\mu_X)$ and $\mathrm{supp}(\mu_Y)$ and $Z$,
  respectively, and a coupling $\hat{\mu}$ of $(\varphi_{X})_\ast\mu_X$
  and $(\varphi_{Y})_\ast\mu_Y$ such that
\be{p6q}
   \hat{\mu}\big\{(z,z'):\,r_Z(z,z')\ge\delta\big\}
 <
   \delta.
\ee

Hence with the special choice of a coupling $\tilde\mu$ of $\mu_X$ and
$\mu_Y$ defined by $\tilde{\mu}(A\times
B)=\hat{\mu}\big(\varphi_{X}(A)\times \varphi_{Y}(B)\big)$, for all
$A\in{\mathcal B}(X)$ and $B\in{\mathcal B}(Y)$,
\be{l15}
\begin {aligned}
   &\tilde\mu^{\otimes 2}\big\{(x,y),(x',y')\in X\times Y:\, |r_X(x,x') -
   r_Y(y,y')|\ge
   2\delta\}
  \\
 &\leq
   \tilde\mu^{\otimes 2}\big\{(x,y), (x',y')\in X\times Y:\,
  \\
 &\qquad
     r_Z(\varphi_X(x), \varphi_Y(y))\ge\delta \text{ or }
     r_Z(\varphi_X(x'), \varphi_Y(y'))\ge\delta\big\}
  \\
 &<
   2\delta.
\end{aligned}
\ee
This implies that $d_{\mathrm{Eur}}\big({\mathcal X},{\mathcal Y}\big)<2\delta$.
\end{proof}\sm

\begin{proof}[Proof of Proposition~\ref{propMet2}]
Observe that by Lemma~\ref{l:dvelta}, $v_\delta(\mathcal X)\overset{\delta\to
  0}{\longrightarrow}0$. So Lemma~\ref{LL} implies the {\em equivalence} of $d_{\rm{GPr}}$ and
  $d_{\rm{Eur}}$ once we have shown that $d_{\rm{Eur}}$ is indeed a
  metric.

The \emph{symmetry} is clear.
If ${\mathcal X}$, ${\mathcal Y}\in\mathbb{M}$ are such that
$d_{\rm{Eur}}(\mathcal X,
\mathcal Y)=0$, by equivalence,
$d_{\rm{GPr}}(\mathcal
X,\mathcal Y)=0$ and hence $\mathcal X=\mathcal Y$.

For the \emph{triangle inequality}, let $\mathcal
X_i=(X_i,r_i,\mu_i)\in\mathbb{M}, i=1,2,3$, be such that
$d_{\rm{Eur}}(\mathcal X_1,\mathcal X_2)<\varepsilon$ and
$d_{\rm{Eur}}(\mathcal X_2,\mathcal X_3)<\delta$ for some
$\varepsilon,\delta>0$. Then there exist couplings $\tilde\mu_{1,2}$
of $\mu_1$ and $\mu_2$ and $\tilde\mu_{2,3}$ of $\mu_2$ and $\mu_3$
with
\be{cou12}
   \tilde\mu_{1,2}^{\otimes 2}\big\{(x_1,x_2),(x_1', x_2'):\,
   |r_1(x_1,x_1') - r_2(x_2,x_2')|\ge\varepsilon\big\}
 <
   \varepsilon
\end{equation}
and
\be{cou23}
   \tilde\mu_{2,3}^{\otimes 2}\big\{(x_2,x_3), (x_2', x_3'):\,
   |r_2(x_2,x_2') - r_3(x_3,x_3')|\ge\delta\big\}
 <
   \delta.
\end{equation}

Introduce the transition kernel $K_{2,3}$ from $X_2$ to
$X_3$ defined by
\be{kernel}
   \tilde\mu_{2,3}({\rm d}(x_2,x_3))
 =
   \mu_2({\rm d}x_2) K_{2,3}(x_2, {\rm d}x_3).
\ee
which exists since $X_2$ and $X_3$ are Polish.

Using this kernel, define a coupling $\tilde\mu_{1,3}$ of $\mu_1$ and
$\mu_3$ by
\begin{equation}
  \tilde\mu_{1,3} ({\rm d}(x_1,x_3))
 :=
  \int_{X_2} \tilde\mu_{1,2}({\rm d}(x_1,x_2)) K_{2,3}(x_2,{\rm d}x_3).
\end{equation}

Then
\begin{equation}
\begin{aligned}
   \tilde\mu_{1,3}^{\otimes 2}&\big\{ (x_1,x_3),(x_1',x_3'):\, |r_1(x_1,x_1')
   - r_3(x_3,x_3')|\ge\varepsilon + \delta\big\}
  \\
 & =
   \int_{X_1^2\times X_2^2\times X_3^2} \tilde\mu_{1,2}({\rm
    d}(x_1,x_2))\tilde\mu_{1,2}({\rm d}(x_1',x_2')) K_{2,3}(x_2,{\rm
    d}x_3) K_{2,3}(x_2',{\rm d}x_3')
  \\
 &\qquad \qquad \qquad \qquad\qquad \qquad \qquad \quad
    \mathbf{1}\big\{|r_1(x_1,x_1') -
  r_3(x_3,x_3')|\ge\varepsilon + \delta\big\}
  \\
 &\leq
   \int_{X_1^2\times
    X_2^2\times X_3^2} \tilde\mu_{1,2}({\rm
    d}(x_1,x_2))\tilde\mu_{1,2}({\rm d}(x_1',x_2')) K_{2,3}(x_2,{\rm
    d}x_3) K_{2,3}(x_2',{\rm d}x_3') \\& \qquad \big(
  \mathbf{1}\{|r_1(x_1,x_1') - r_2(x_2,x_2')|\ge\varepsilon\} + \mathbf{1}\{|r_2(x_2,x_2')
  - r_3(x_3,x_3')|\ge\delta\}\big)
   \\
 &=
   \tilde\mu_{1,2}^{\otimes 2}\big\{(x_1,x_2), (x_1 ,x_2'):\, |r_1(x_1,x_1') - r_2(x_2,x_2')| \ge
  \varepsilon\} \\&\qquad \qquad \qquad + \tilde\mu_{2,3}^{\otimes 2}
  \{(x_2,x_3), (x_2 ,x_3'):\, |r_2(x_2,x_2') - r_3(x_3,x_3')|\ge \delta\big\}
  \\
 &<
   \varepsilon + \delta
\end{aligned}
\end{equation}
which yields $d_{\rm{Eur}}(\mathcal X_1, \mathcal X_3)<\varepsilon+\delta$.
\end{proof}\sm

\subsection*{The Gromov-Wasserstein and the modified Eurandom metric}
The topology of weak convergence for probability measures on a fixed
metric space $(Z,r)$ is generated not only by the Prohorov metric, but
also by
\be{Was}
   d^{(Z,r_Z)}_{\mathrm{W}}(\mu_1,\mu_2)
 :=
   \inf_{\tilde\mu}\int_{Z\times Z}\tilde\mu({\rm d}(x,x'))\,\big(r(x,x')\wedge 1\big),
\end{equation}
where the infimum is over all couplings $\tilde\mu$ of $\mu_1$ and
$\mu_2$. This is a version of the {\em Wasserstein metric} (see, for example,
\cite{Rachev1991}). If we rely on the Wasserstein rather than the
Prohorov metric, this results in two further metrics:
in the Gromov-Wasserstein metric,
i.e.,
\be{GW}
   d_{\rm{GW}}(\mathcal X, \mathcal Y)
 :=
   \inf_{(\varphi_X, \varphi_Y,Z)}d_{\rm{W}}^{(Z,r_Z)}
   \big((\varphi_X)_\ast\mu_X,(\varphi_Y)_\ast\mu_Y\big),
\end{equation}
where the infimum is over all isometric embeddings from $\mathrm{supp}(\mu_X)$ and $\mathrm{supp}(\mu_Y)$ into a common
metric $Z$ and in the modified Eurandom metric
\begin{equation}\label{eq:p003}
  \begin{aligned}
    d_{\mathrm{Eur}}'\big(&\mathcal X, \mathcal Y\big):= \\ &
    \inf_{\tilde\mu} \int\tilde\mu(\mathrm{d}(x,y)) \tilde
    \mu(\mathrm{d}(x',y'))\big(|r_X(x,x')-r_Y(y,y')|\wedge 1\big),
  \end{aligned}
\end{equation}
where the infimum is over all couplings of $\mu_X$ and $\mu_Y$.

\begin{remark}
An $L^2$-version of $d_{\rm{GW}}$
on the set of compact metric measure spaces is already used in
\cite{Stu2006}.
It turned out that the
metric is complete and the generated topology is separable.
\hfill$\qed$\end{remark}\sm

Altogether, we might ask if we could achieve similar bounds to those
given in Lemma~\ref{LL}  by exchanging the Gromov-Prohorov with the
Gromov-Wasserstein metric and the Eurandom with the modified
Eurandom metric.
\begin{proposition}
  The distances $d_{\rm{GW}}$ and $d_{\rm{Eur}}'$ define metrics on
  $\mathbb{M}$. They all generate the Gromov-Prohorov topology. \label{P:PP}
  Bounds that relate these two metrics with $d_{\rm{GPr}}$ and
  $d_{\rm{Eur}}$ are for $\mathcal X, \mathcal Y\in\mathbb M$,
  \begin{equation}
    (d_{\rm{GPr}}(\mathcal X, \mathcal Y))^2 \leq d_{\rm{GW}}(\mathcal X, \mathcal Y)
    \leq d_{\rm{GPr}}(\mathcal X, \mathcal Y)
  \end{equation}
  and
  \begin{equation}
    (d_{\rm{Eur}}(\mathcal X, \mathcal Y))^2 \leq
    d_{\rm{Eur}}'(\mathcal X, \mathcal Y) \leq d_{\rm{Eur}}(\mathcal
    X, \mathcal Y)
  \end{equation}

Consequently, the Gromov-Wasserstein metric is complete.
\end{proposition}\sm

\begin{proof}
  The fact that $d_{\rm{GW}}$ and $d_{\rm{Eur}}'$ define metrics on
  $\mathbb{M}$ is proved analogously as for the Gromov-Prohorov and
  the Eurandom metric. The Prohorov and the version of
  the Wasserstein metric used in \eqref{GW} and \eqref{eq:p003} on
  fixed metric spaces can be bounded uniformly (see, for example,
  Theorem~3 in
   \cite{GibbsSu2001}). This immediately carries over to the
  present case.
\end{proof}\sm

\appendix
\section{Additional facts on Gromov-Hausdorff convergence}
Recall the notion of the Gromov-Hausdorff distance on the space $\mathbb X_c$ of isometry classes of compact metric spaces given in \eqref{eq:GH1}. We give a statement concerning convergence in the Gromov-Hausdorff metric which is analogous to Lemma \ref{l:Gpronespace} for Gromov-Prohorov convergence.

\begin{lemma}\label{l:met2}
  Let $(X,r_X)$, $(X_1,r_{X_1})$, $(X_2,r_{X_2})$, ... be in $\mathbb
  X_{\rm{c}}$. Then
  $$d_{\mathrm{GH}}(X_n,X)\overset{n\to\infty}{\longrightarrow} 0$$ if
  and only if there is a compact metric space $(Z,r_Z)$ and isometric
  embeddings $\varphi$, $\varphi_1$, $\varphi_2$, ... of $(X,r)$, $(X_1,r_{X_1})$,
  $(X_2,r_{X_2})$, ..., respectively, into $(Z,r_Z)$ such that
  \begin{align}\label{eq:A21}d^{(Z,r_Z)}_{\mathrm{H}}\big(\varphi_n(X_n),\varphi(X)\big)
 \overset{n\to\infty}{\longrightarrow}
  0.
  \end{align}
\end{lemma}\sm

\begin{proof}The ``if''-direction is clear. So we come immediately to
  the ``only if'' direction. If
  $d_{\mathrm{GH}}\big(X_n,X\big)\overset{n\to\infty}{\longrightarrow}
  0$, then by (\ref{GH}) we find correspondences $R_n$ between $X$
  and $X_n$ such that
  $\mathrm{dis}(R_n)\overset{n\to\infty}{\longrightarrow} 0$.  Using these and $X_0:=X$, we define recursively metrics $r_{Z_n}$ on $Z_n:=
  \bigsqcup\nolimits_{k=0}^n X_k$. First, set $Z_1:=X_0\sqcup X_1$ and
  $r_{Z_1} := r_{Z_1}^{R_1}$ (recall Remark \ref{rem:met1}). In the $n^{\mathrm{th}}$ step, we are
  given a metric on $Z_n$. Consider the canonical isometric embedding $\varphi$ from $X$ to
  $Z_n$ and define the relation $\tilde R_{n}\subseteq Z_n\times
  X_{n+1}$ by
\be{Rn1}
   \tilde R_{n+1}
 :=
   \big\{(z,x_{})\in Z_n\times X_{n+1}:\,(\varphi^{-1}(z),x_{}) \in R_{n+1}\big\},
\end{equation}
and set $r_{Z_{n+1}}:=r_{Z_{n+1}}^{\tilde R_{n+1}}$. By this procedure
we end up with a metric $r_{Z}$ on $Z:=\bigsqcup\nolimits_{n=0}^\infty X_n$ and isometric embeddings $\varphi_0$,
$\varphi_1$,$...$ between $X_0$, $X_1$, ... and
$Z$, respectively, such that
\be{e:009a}
   d^{(Z,r_Z)}_{\mathrm{H}}\big(\varphi_n(X_n),\varphi(X)\big)
 =
   \frac12 \mathrm{dis}(R_n) \overset{n\to\infty}{\longrightarrow} 0.
\end{equation}
W.l.o.g.\ we can assume that $Z$ is complete. Otherwise we just embed
everything into the completion of $Z$. To verify {\em compactness} of $(Z,r_Z)$
it is therefore sufficient to
show that $Z$ is totally bounded (see, for example, Theorem 1.6.5 in
\cite{BurBurIva01}).  For that purpose fix $\varepsilon>0$, and let
$n\in\N$. Since $X$ is compact, we can choose a finite
$\varepsilon/2$-net $S$ in $X$. Then for all $x\in Z$ with
$r_Z(x,X)<\varepsilon/2$ there exists $x'\in S$ such that
$r_Z(x,x')<\varepsilon$. Moreover,
$d_{\mathrm{H}}\big(\varphi_n(X_n),\varphi(X)\big)<\varepsilon$, for
all but finitely many $n\in\N$.
For the remaining $\varphi_n(X_n)$
choose finite $\varepsilon$-nets and denote their union by $\tilde S$.
In this way, $S\cup \tilde S$ is a finite set, and
$\{B_\varepsilon(s):\,s\in S\cup \tilde S\}$ is a covering of $Z$.
\end{proof}\sm

{\sc Acknowledgements. } The authors thank Anja Sturm and Theo Sturm
for helpful discussions, and Reinhard Leipert for help on
Remark~\ref{Rem:07}(ii). Our special thanks go to Steve Evans for
suggesting to verify that the Gromov-weak topology is not weaker than
the Gromov-Prohorov topology and to Vlada Limic for encouraging us to
write a paper solely on the topological aspects of genealogies.
Finally we thank an anonymous referee for several helpful comments
that improved the presentation and correctness of the paper.

\newcommand{\etalchar}[1]{$^{#1}$}


\end{document}